\documentclass[final,3p]{elsarticle}
 \usepackage[latin1]{inputenc}
 \usepackage{graphics}
 \usepackage{graphicx}
 \usepackage{epsfig}
\usepackage{amssymb}
 \usepackage{amsthm}
 \usepackage{lineno}
 \usepackage{amsmath}
   \numberwithin{equation}{section}
\usepackage{mathrsfs}
\usepackage{esint}

\newtheorem{thm}{Theorem}[section]

\newtheorem{lem}[thm]{Lemma}

\newtheorem{defn}[thm]{Definition}

\newtheorem{exam}[thm]{Example}

\biboptions{sort&compress,square}
\begin{document}
\begin{frontmatter}
\author[rvt1]{Jian Wang}
\ead{wangj@tute.edu.cn}
\author[rvt2]{Yong Wang\corref{cor2}}
\ead{wangy581@nenu.edu.cn}
\author[rvt1]{Mingyu Liu}
\ead{13246450531@163.com}
\cortext[cor2]{Corresponding author.}
\address[rvt1]{School of Science, Tianjin University of Technology and Education, Tianjin, 300222, P.R.China}
\address[rvt2]{School of Mathematics and Statistics, Northeast Normal University,
Changchun, 130024, P.R.China}

\title{Spectral forms and de-Rham Hodge operator}
\begin{abstract}
Motivated by the trilinear functional of differential one-forms, spectral triple and spectral torsion
 for the Hodge-Dirac operator, we introduce a multilinear functional of differential one-forms for a
finitely summable regular spectral triple with a noncommutative residue,
which generalize the spectral torsion defined by  Dabrowski-Sitarz-Zalecki.
The main results of this paper recover two forms, torsion of the linear connection and four forms by the noncommutative residue
 and perturbed de-Rham Hodge operators,
and provides an explicit computation of
generalized spectral forms associated with the perturbed de-Rham Hodge Dirac triple.
\end{abstract}
\end{frontmatter}
\section{Introduction}
\label{1}

In Connes' program of noncommutative geometry, the role of geometrical
objects is played by spectral triples $(\mathcal{A},\mathcal{H},D)$. Similar to the commutative case and
the canonical spectral triple $(C^{\infty}(M),L^{2}(S),D)$, where $(M,g,S)$ is a closed spin
manifold and $D$ is the Dirac operator acting on the spinor bundle $S$, the spectrum
of the Dirac operator $D$ of a spectral triple $(\mathcal{A},\mathcal{H},D)$ encodes the geometrical
information of the spectral triple.
The noncommutative residue of Wodzicki \cite{Wo,Wo1} and Guillemin \cite{VG} is a trace
on the algebra of (integer order) $\Psi DOs$ on $M$. An important feature is that it
allows us to extend to all $\Psi DOs$ the Dixmier trace, which plays the role of the
integral in the framework of noncommutative geometry.
This extends the Wodzicki residue from pseudodif-
ferential operators on a manifold to the general framework of spectral triples, and
gives meaning to $ \fint T$ for the compact operators $T$. It is simply given by
\begin{align}
\fint T=Res_{s=0} {\rm{Tr}} (T|D|^{-s}).
\end{align}
 In\cite{MA,Wo,Wo1, Ac} it was shown that the Wodzicki
residue ${\rm  Wres}(\Delta^{-n/2+1})$ of a generalized laplacian $\Delta$  on a complex vector bundle $E$ over a closed compact manifold $M$,
is the integral of the second coefficient of the heat kernel expansion of $\Delta$ up to a proportional factor.
In \cite{Co1, Co2,Co,CCS}, Connes used the noncommutative residue to derive a conformal 4-dimensional Polyakov action analogue,
and made a challenging observation that the noncommutative residue of the square of the inverse of the
Dirac operator was proportional to the Einstein-Hilbert action. Then the Kastler-Kalau-Walze type
theorem \cite{KW,Ka} gives an operator-theoretic explanation of the gravitational action for a  lower dimensional closed spin manifold.
 Furthermore, Fedosov et al. defined the noncommutative residue on Boutet de Monvel's algebra and proved that it was a
unique continuous trace for the case of manifolds with boundary \cite{FGLS}. In \cite{S}, Schrohe gave the relation between the Dixmier trace
and the noncommutative residue for manifolds with boundary.
Wang generalized Connes' results to the case of manifolds with boundary,
and proved  Kastler-Kalau-Walze type theorems for lower-dimensional manifolds
with boundary associated with Dirac operator and signature operator \cite{Wa1,Wa2,Wa3}.

The trilinear functional of differential one-forms for a finitely
summable regular spectral triple with a noncommutative residue has been computed recently in \cite{DSZ},
Dabrowski et al. demonstrated that for
a canonical spectral triple over a closed spin manifold it recovers the torsion of the linear
connection, which is a first step towards linking the spectral approach with the algebraic approach based
on Levi-Civita connections.
 It seems also natural to extend this approach to the metric connection and bring
into focus a non-zero torsion.
The importance of the spectral torsion has been stressed recently in a somewhat different
context. For example, Dirac operators with torsion are by now well-established analytical tools in the study of special geometric structures.
 Ackermann and  Tolksdorf \cite{AT} proved a generalized version of the well-known Lichnerowicz formula for the square of the
most general Dirac operator with torsion  $D_{T}$ on an even-dimensional spin manifold associated to a metric connection with torsion.
In \cite{PS,PS1}, Pf$\ddot{a}$ffle and Stephan considered orthogonal connections with arbitrary torsion on compact Riemannian manifolds,
and for the induced Dirac operators, twisted Dirac operators and Dirac operators of Chamseddine-Connes type they computed the spectral
action.
 Sitarz and Zajac \cite{SZ} investigated the spectral action for scalar perturbations of Dirac operators.
Iochum and Levy \cite{IL}computed the heat kernel coefficients for Dirac operators with one-form perturbations.
Wang \cite{Wa4} considered the arbitrary perturbations of Dirac operators, and established the associated Kastler-
Kalau-Walze theorem.
In \cite{WWW,WW2,WW3}, we gave two Lichnerowicz type formulas for Dirac operators
and signature operators twisted by a vector bundle with a non-unitary connection,  and proved the Kastler-Kalau-Walze type theorem associated to
Dirac operators with torsion on compact manifolds with boundary.
In \cite{LWW1}, Liu et al.  proved the Kastler-Kalau-Walze type theorems for the perturbation of the de Rham Hodge operator on 4-dimensional
and 6-dimensional compact manifolds.

Based on the spectral torsion and the noncommutative residue, Dabrowski et al. \cite{DSZ} showed that the spectral definition of
 torsion can be readily extended to the noncommutative case of spectral triples.
By twisting the spectral triple of a Riemannian spin manifold, Martinetti et al. showed how
to generate an orthogonal and geodesic preserving torsion from a torsionless Dirac operator in \cite{MNZ}.
Motivated by the spectral torsion  \cite{DSZ}, the de-Rham Hodge operators with torsion and the noncommutative residue,
it is a natural question to define the spectral forms, which generalize the spectral torsion.
We introduce a multilinear functional of differential one-forms for a
finitely summable regular spectral triple with a noncommutative residue,
which generalize the spectral torsion defined by  Dabrowski-Sitarz-Zalecki.
We reconsider the calculations of the spectral torsion as three forms  associated with the de-Rham Hodge operators with torsion.
The purpose of the present paper is to  recover two forms, torsion of the linear connection and four forms,
and provide an explicit computation of
generalized spectral torsion associated with the perturbed de-Rham Hodge Dirac triple.
 The aim of this paper is to prove the following Theorems.
\begin{thm}
For the de-Rham Hodge operator with the bilinear Clifford multiplication by functional of differential one-forms
$\widehat{c}(u),\widehat{c}(v)$, the spectral torsion for  $ d+\delta+\sqrt{-1}T_{1} $  equals to
\begin{align}
 \mathscr{T}_{1}(\widehat{c}(u),\widehat{c}(v))
=&\int_M\int_{|\xi|=1}{\rm Tr}_{\wedge^{*}(T^{*}M)}
[\sigma_{-2m}\big(\widehat{c}(u)\widehat{c}(v) (d+\delta+\sqrt{-1}T_{1}) ^{-2m+1}\big)]\sigma(\xi){\rm d}x\nonumber\\
=&\int_M (2m-1)2^{2m} \sqrt{-1}T(u,v )
\rm{vol}(S^{2m-1}) {\rm d}\rm{vol}_{M},
\end{align}
where $T_{1}=\sum_{k<l }T(\widetilde{e_k},\widetilde{e_l})
\widehat{c}(\widetilde{e_j})\widehat{c}(\widetilde{e_l})$ is a two form.
\end{thm}
\begin{thm}
For the de-Rham Hodge operator with the trilinear Clifford multiplication by functional of differential one-forms
$c(u),c(v),c(w)$, the spectral torsion for  $ d+\delta+T_{2} $   equals to
\begin{align}
&\mathscr{T}_{2}(c(u),c(v),c(w))\nonumber\\
=&\int_M\int_{|\xi|=1}{\rm Tr}_{\wedge^{*}(T^{*}M)}
[\sigma_{-2m}\big(c(u)c(v)c(w)(d+\delta+T_{2})^{-2m+1}\big)]\sigma(\xi){\rm d}x\nonumber\\
=&\int_M (3-18m)2^{2m-1} T(u,v,w)
\rm{vol}(S^{2m-1}) {\rm d}\rm{vol}_{M},
\end{align}
where $T_{2}=\frac{3}{2}\sum_{i<s<t}T(\widetilde{e_i},\widetilde{e_s},\widetilde{e_t})
c(\widetilde{e_i})c(\widetilde{e_s})c(\widetilde{e_t})
-\frac{1}{4}\sum^n_{i\neq s\neq t}T(\widetilde{e_i},\widetilde{e_s},\widetilde{e_t})
c(\widetilde{e_i})\widehat{c}(\widetilde{e_s})\widehat{c}(\widetilde{e_t})$ denotes the `lifted' torsion term.
\end{thm}
\begin{thm}
For the de-Rham Hodge operator with the trilinear Clifford multiplication by functional of differential one-forms
$c(u),\widehat{c}(v),\widehat{c}(w)$, the spectral torsion for  $ d+\delta+T_{2} $   equals to
\begin{align}
&\mathscr{T}_{3}(c(u),\widehat{c}(v),\widehat{c}(w))\nonumber\\
=&\int_M\int_{|\xi|=1}{\rm Tr}_{\wedge^{*}(T^{*}M)}
[\sigma_{-2m}\big(c(u)\widehat{c}(v)\widehat{c}(w)(d+\delta+T_{2})^{-2m+1}\big)]\sigma(\xi){\rm d}x\nonumber\\
=&\int_M (1-2m)2^{2m-1}T(u,v,w)
\rm{vol}(S^{2m-1}) {\rm d}\rm{vol}_{M},
\end{align}
where $T_{2}=\frac{3}{2}\sum_{i<s<t}T(\widetilde{e_i},\widetilde{e_s},\widetilde{e_t})
c(\widetilde{e_i})c(\widetilde{e_s})c(\widetilde{e_t})
-\frac{1}{4}\sum^n_{i\neq s\neq t}T(\widetilde{e_i},\widetilde{e_s},\widetilde{e_t})
c(\widetilde{e_i})\widehat{c}(\widetilde{e_s})\widehat{c}(\widetilde{e_t})$ denotes the `lifted' torsion term..
\end{thm}
By Theorem 1.2 and Theorem 1.3, we can extend the spectral torsion for de-Rham Hodge operators to
de-Rham Hodge spectral triple with four-form perturbations.

\begin{thm}
For the de-Rham Hodge operator with the four-linear Clifford multiplication by functional of differential one-forms
$c(u),c(v),\widehat{c}(w),\widehat{c}(z)$, the spectral torsion for  $ d+\delta+\sqrt{-1}T_{3} $   equals to
\begin{align}
&\mathscr{T}_{4}(c(u),c(v),\widehat{c}(w) ,\widehat{c}(z))\nonumber\\
=&\int_M\int_{|\xi|=1}{\rm Tr}_{\wedge^{*}(T^{*}M)}
[\sigma_{-2m}\big(c(u)c(v)\widehat{c}(w)\widehat{c}(z)(d+\delta+\sqrt{-1}T_{3} ) ^{-2m+1}\big)]\sigma(\xi){\rm d}x\nonumber\\
=&\int_M \frac{(-2m-1)}{3}2^{2m-1}\sqrt{-1} T(u,v,w,z)
\rm{vol}(S^{n-1}) {\rm d}\rm{vol}_{M}.
\end{align}
where $T_{3}=\sum_{k<l <\alpha<\beta}T(\widetilde{e_k},\widetilde{e_l},\widetilde{e_\alpha},\widetilde{e_\beta})
c(\widetilde{e_k})c(\widetilde{e_l})\widehat{c}(\widetilde{e_\alpha})\widehat{c}(\widetilde{e_\beta})$ is a four form.
\end{thm}

\begin{thm}
For the de-Rham Hodge operator with the  four-linear Clifford multiplication by functional of differential one-forms
$\widehat{c}(u),\widehat{c}(v),\widehat{c}(w),\widehat{c}(z)$, the spectral torsion for  $ d+\delta+\sqrt{-1}T_{4} $   equals to
\begin{align}
&\mathscr{T}_{5}(\widehat{c}(u) ,\widehat{c}(v),\widehat{c}(w),\widehat{c}(z))\nonumber\\
=&\int_M\int_{|\xi|=1}{\rm Tr}_{\wedge^{*}(T^{*}M)}
[\sigma_{-2m}\big(\widehat{c}(u)\widehat{c}(v)\widehat{c}(w)\widehat{c}(z)(d+\delta+\sqrt{-1}T_{4})^{-2m+1}\big)]\sigma(\xi){\rm d}x\nonumber\\
=&\int_M (-2m+1)2^{2m}\sqrt{-1} T(u,v,w,z)
\rm{vol}(S^{n-1}) {\rm d}\rm{vol}_{M}.
\end{align}
where $T_{4}=\sum_{k<l <\alpha<\beta}T(\widetilde{e_k},\widetilde{e_l},\widetilde{e_\alpha},\widetilde{e_\beta})
\widehat{c}(\widetilde{e_k})\widehat{c}(\widetilde{e_l})\widehat{c}(\widetilde{e_\alpha})\widehat{c}(\widetilde{e_\beta})$ is a four form.
\end{thm}

\section{ spectral two-forms for the de-Rham triple}

The purpose of this section is to specify the de-Rham Hodge type operator and compute the spectral two-forms for the de-Rham triple.
Let $M$ be a $n$-dimensional ($n\geq 3$) oriented compact Riemannian manifold with a Riemannian metric $g^{M}$.
 And let $\nabla^L$ be the Levi-Civita connection about $g^M$. In the local coordinates $\{x_i; 1\leq i\leq n\}$ and the
fixed orthonormal frame $\{\widetilde{e_1},\cdots,\widetilde{e_n}\}$, the connection matrix $(\omega_{s,t})$ is defined by
\begin{equation}
\nabla^L(\widetilde{e_1},\cdots,\widetilde{e_n})= (\widetilde{e_1},\cdots,\widetilde{e_n})(\omega_{s,t}).
\end{equation}
\indent Let $\epsilon (\widetilde{e_j*})$,~$\iota (\widetilde{e_j*})$ be the exterior and interior multiplications respectively, $\widetilde{e_j*}$ be the dual base of $\widetilde{e_j}$  and $c(\widetilde{e_j})$ be the Clifford action.
Suppose that $\partial_{i}$ is a natural local frame on $TM$
and $(g^{ij})_{1\leq i,j\leq n}$ is the inverse matrix associated to the metric
matrix  $(g_{ij})_{1\leq i,j\leq n}$ on $M$. Write
\begin{equation}
c(\widetilde{e_j})=\epsilon (\widetilde{e_j*})-\iota (\widetilde{e_j*});\
\widehat{c}(\widetilde{e_j})=\epsilon (\widetilde{e_j*})+\iota (\widetilde{e_j*}),
\end{equation}
which satisfies
\begin{align}
&c(\widetilde{e_i})c(\widetilde{e_j})+c(\widetilde{e_j})c(\widetilde{e_i})=-2\delta_i^j;\nonumber\\
&\widehat{c}(\widetilde{e_i})c(\widetilde{e_j})+c(\widetilde{e_j})\widehat{c}(\widetilde{e_i})=0;\nonumber\\
&\widehat{c}(\widetilde{e_i})\widehat{c}(\widetilde{e_j})+\widehat{c}(\widetilde{e_j})\widehat{c}(\widetilde{e_i})=2\delta_i^j.
\end{align}
By \cite{Wa3}, we have the signature operator
\begin{align}
 d+\delta =\sum^n_{i=1}c(\widetilde{e_i})\bigg[\widetilde{e_i}-\frac{1}{4}
\sum_{s,t}\langle \nabla^{L}_{\widetilde{e_i}}\widetilde{e_s}, \widetilde{e_t} \rangle
\Big[\widehat{c}(\widetilde{e_s})\widehat{c}(\widetilde{e_t})-c(\widetilde{e_s})c(\widetilde{e_t})\Big]\bigg].
\end{align}
Let $g^{ij}=g(dx_{i},dx_{j})$, $\xi=\sum_{k}\xi_{j}dx_{j}$ and $\nabla^L_{\partial_{i}}\partial_{j}=\sum_{k}\Gamma_{ij}^{k}\partial_{k}$.
We collect here, for the reader's convenience, all necessary facts on the explicit representation of $(d+\delta+A)^{2}$:
\begin{align}
&\sigma_{i}=-\frac{1}{4}\sum_{s,t}\omega_{s,t}(\widetilde{e_i})c(\widetilde{e_s})c(\widetilde{e_t});
\ a_{i}=\frac{1}{4}\sum_{s,t}\omega_{s,t}(\widetilde{e_i})\widehat{c}(\widetilde{e_s})\widehat{c}(\widetilde{e_t});\nonumber\\
&\xi^{j}=g^{ij}\xi_{i};\ \Gamma^{k}=g^{ij}\Gamma_{ij}^{k};\ \sigma^{j}=g^{ij}\sigma_{i};\ a^{j}=g^{ij}a_{i};\nonumber\\
&A=\sum_{i_{1}<i_{2}<\cdots<i_{n}}T(\widetilde{e_{i_{1}}},\widetilde{e_{i_{2}}},\cdots,\widetilde{e_{i_{n}}})
c(\widetilde{e_{i_{1}}})c(\widetilde{e_{i_{2}}})\cdots c(\widetilde{e_{i_{n}}}).
\end{align}
Then, the perturbation of the de Rham Hodge operator $D_{A}$ can be written as
\begin{align}
{D}_{A}&=\sum^n_{i=1}c(\widetilde{e_i})[\widetilde{e_i}+\sigma_{i}+a_{i}]+A.
\end{align}
From \cite{Ac} and \cite{Y}, the local expression of $(d+\delta)^2$ is
\begin{align}
(d+\delta)^2=-\Delta_0-\frac{1}{8}\sum_{i,j,k,l}R_{ijkl}\widehat{c}(\widetilde{e_i})
\widehat{c}(\widetilde{e_j})c(\widetilde{e_k})c(\widetilde{e_l})+\frac{1}{4}s,
\end{align}
where
\begin{align}
-\Delta_0=-g^{ij}(\nabla_i^L\nabla_j^L-\Gamma_{ij}^k\nabla_k^L).
\end{align}
An easy calculation gives
\begin{align}
{D}^2_{A}&=(d+\delta)^2+(d+\delta)A+A(d+\delta)+A^2.
\end{align}
Also, straightforward computations yield
\begin{align}
(d+\delta)A+A(d+\delta)=&\sum_{i,j}g^{ij}[c(\partial_{i})A+Ac(\partial_{i})]\partial_{j}+\sum_{i,j}g^{ij}[c(\partial_{i})\partial_{j}(A)\nonumber\\
&+c(\partial_{i})a_{j}A+c(\partial_{i})\sigma_{j}A+Ac(\partial_{i})a_{j}+Ac(\partial_{i})\sigma_{j}].
\end{align}
In what follows, using the Einstein sum convention for repeated index summation:
$\partial^j=g^{ij}\partial_i$,$\sigma^i=g^{ij}\sigma_{j}$, $\Gamma_{k}=g^{ij}\Gamma_{ij}^{k}$,
 and we omit the summation sign, from (6a) in \cite{Ka} and (2.23) in \cite{LWW1}, we have
\begin{align}
{D}^2_{A}&=-\sum_{i,j}g^{ij}[\partial_{i}\partial_{j}+2\sigma_{i}\partial_{j}+2a_{i}\partial_{j}-\Gamma_{ij}^{k}\partial_{k}
+(\partial_{i}\sigma_{j})+(\partial_{i}a_{j})+\sigma_{i}\sigma_{j}+\sigma_{i}a_{j}+a_{i}\sigma_{j}\nonumber\\
&+a_{i}a_{j}-\Gamma_{ij}^{k}\sigma_{k}-\Gamma_{ij}^{k}a_{k}]+\sum_{i,j}g^{ij}[c(\partial_{i})\partial_{j}(A)
+c(\partial_{i})a_{j}A+c(\partial_{i})\sigma_{j}A+Ac(\partial_{i})a_{j}\nonumber\\
&+Ac(\partial_{i})\sigma_{j}]+\sum_{i,j}g^{ij}[c(\partial_{i})A+Ac(\partial_{i})]\partial_{j}
-\frac{1}{8}\sum_{i,j,k,l}R_{ijkl}\widehat{c}(\widetilde{e_i})\widehat{c}(\widetilde{e_j})c(\widetilde{e_k})c(\widetilde{e_l})\nonumber\\
&+\frac{1}{4}s+A^2.
\end{align}

Using an explicit formula for $d+\delta$, we can reformulate $d+\delta +T$ as follows.
\begin{defn}
The de-Rham Hodge type operator  is the first order differential operator on $ \Gamma(M,\wedge^{*}(T^{*}M))$  given by the formula
 \begin{align}
d+\delta +T=&d+\delta
+\sum_{i_{1}<i_{2}<\cdots<i_{l}}T(\widetilde{e_{i_{1}}},\widetilde{e_{i_{2}}},\cdots,\widetilde{e_{i_{l}}})
c(\widetilde{e_{i_{1}}})c(\widetilde{e_{i_{2}}})\cdots c(\widetilde{e_{i_{l}}}),
\end{align}
where $T$ is a degree $l-$form.
\end{defn}

In this section we want to consider the bilinear functional  for the  de-Rham Hodge type operator.
To avoid technique terminology we only state our results for
compact oriented Riemannian manifold of even dimension $n=2m$ by using the trace of  de-Rham Hodge type operator and the noncommutative residue density.
\begin{defn}\cite{DSZ}
For $d+\delta +\sqrt{-1}T_{1}$ given by $T_{1}=\sum_{k<l }T(\widetilde{e_k},\widetilde{e_l})
\widehat{c}(\widetilde{e_j})\widehat{c}(\widetilde{e_l})$, the trilinear Clifford multiplication by functional of differential one-forms
$\widehat{c}(u),\widehat{c}(v)$
\begin{align}
\mathscr{T}_{1}(\widehat{c}(u),\widehat{c}(v))={\rm Wres}\big(\widehat{c}(u)\widehat{c}(v)(d+\delta+\sqrt{-1}T_{1})^{-n+1}\big)
\end{align}
is called the spectral two form $\mathscr{T}_{1}$ associated with the de-Rham Hodge type operator.
\end{defn}
The following noncommutative residue formula \cite{DSZ} plays a key role in our proof
of the spectral torsion  associated with the de-Rham Hodge type operator.
\begin{align}
{\rm Wres}\big(\widehat{c}(u)\widehat{c}(v)(d+\delta+\sqrt{-1}T_{1})^{-n+1}\big)=\int_M\int_{|\xi|=1}{\rm Tr}_{\wedge^{*}(T^{*}M)}
[\sigma_{-2m}\big(\widehat{c}(u)\widehat{c}(v)(d+\delta+\sqrt{-1}T_{1})^{-2m+1}\big)]\sigma(\xi){\rm d}x.
\end{align}
Now we consider the explicit representation of the symbols of $\sigma_{-2m}\big(\widehat{c}(u)\widehat{c}(v)(d+\delta+\sqrt{-1}T_{1})^{-2m+1})$.
Let the cotangent vector $\xi=\sum \xi_jdx_j$ and
$\xi^j=g^{ij}\xi_i$.
Decompose the operators $\widetilde{D}_{T}^2=(d+\delta+\sqrt{-1}T_{1}) ^{2}$ and $\widetilde{D}_{T}=(d+\delta+\sqrt{-1}T_{1}) $ by different orders as
\begin{align}
\sigma(\widetilde{D}_{T}^{2})=\sigma_{2}(\widetilde{D}_{T}^{2})
+\sigma_{1}(\widetilde{D}_{T}^{2})+\sigma_{0}(\widetilde{D}_{T}^{2});~~
 \sigma(\widetilde{D}_{T}^{1})=\sigma_{1}(\widetilde{D}_{T}^{1})+\sigma_{0}(\widetilde{D}_{T}^{1}).
\end{align}
Note that the arguments in the proof of Lemma 3.1 in \cite{WW1} can be used to compute
the symbolic representation for the de-Rham Hodge operator without any change, where the same conclusion was derived for de-Rham Hodge operator.
By Lemma 3.1 \cite{WW1}, we get
 \begin{align}
\sigma_{-2}(\widetilde{D}_{T}^{-2})=|\xi|^{-2},~~\sigma_{-(2m-2)}(\widetilde{D}_{T}^{-(2m-2)})=(|\xi|^2)^{1-m},
~~\sigma_{-2m+1}(\widetilde{D}_{T}^{-2m+1})=\frac{\sqrt{-1}c(\xi)}{|\xi|^{2m}}.
\end{align}
Then by (3.8) in \cite{WW1}, we have
 \begin{align}
 \sigma_{-1-2m}(\widetilde{D}_{T}^{-2m})=&m\sigma_2(\widetilde{D}_{T}^2)^{(-m+1)}\sigma_{-3}(\widetilde{D}_{T}^{-2}) \nonumber\\
 &-\sqrt{-1} \sum_{k=0}^{m-2}\sum_{\mu=1}^{2m+2}
\partial_{\xi_{\mu}}\sigma_{2}^{-m+k+1}(\widetilde{D}_{T}^2)
\partial_{x_{\mu}}\sigma_{2}^{-1}(\widetilde{D}_{T}^2)(\sigma_2(\widetilde{D}_{T}^2))^{-k}.
\end{align}
In the same way we obtain
\begin{align}
 \sigma_{-2m}(\widetilde{D}_{T}^{1-2m})=&\sigma_{-2m}(\widetilde{D}_{T}^{-2m}\cdot \widetilde{D}_{T})
=\Big\{\sum_{|\alpha|=0}^{+\infty}(-\sqrt{-1})^{|\alpha|}
\frac{1}{\alpha!}\partial^\alpha_\xi[\sigma(\widetilde{D}_{T}^{-2m})]\partial^\alpha_x[\sigma(\widetilde{D}_{T})]\Big\}_{-2m}   \nonumber\\
 =&\sigma_{-2m}(\widetilde{D}_{T}^{-2m})\sigma_0(\widetilde{D}_{T})+\sigma_{-2m-1}(\widetilde{D}_{T}^{-2m})\sigma_1(\widetilde{D}_{T})\nonumber\\
&+\sum_{|\alpha|=1}(-\sqrt{-1})
\partial^\alpha_\xi[\sigma_{-2m}(\widetilde{D}_{T}^{-2m})]\partial^\alpha_x[\sigma_1(\widetilde{D}_{T})]   \nonumber\\
 =&|\xi|^{-2m}\sigma_0(\widetilde{D}_{T})+\sum_{j=1}^{2m+2}\partial_{\xi_j}(|\xi|^{-2m})\partial_{x_j}c(\xi)+
\Big[m\sigma_2(\widetilde{D}_{T}^2)^{-m+1}\sigma_{-3}(\widetilde{D}_{T}^{-2})  \nonumber\\
&-\sqrt{-1} \sum_{k=0}^{m-2}\sum_{\mu=1}^{2m+2}
\partial_{\xi_{\mu}}\sigma_{2}^{-m+k+1}(\widetilde{D}_{T}^2)
\partial_{x_{\mu}}\sigma_{2}^{-1}(\widetilde{D}_{T}^2)(\sigma_2(\widetilde{D}_{T}^2))^{-k}\Big ]\sqrt{-1}c(\xi).
\end{align}
For any fixed point $x_0\in M$, we can choose the normal coordinates
$U$ of $x_0$ in $M$ and compute $\sigma_{-n}\big(\widehat{c}(u)\widehat{c}(v) (d+\delta+\sqrt{-1}T_{1})^{-n+1}\big)$.
Let $c(\widetilde{e_j})=\epsilon (\widetilde{e_j*})-\iota (\widetilde{e_j*});\
\widehat{c}(\widetilde{e_j})=\epsilon (\widetilde{e_j*})+\iota (\widetilde{e_j*})$ act on $\Gamma(M,\wedge^{*}(T^{*}M))$.
Then we have $ \partial_{x_j}[c(\widetilde{e_{i}})]=\partial_{x_j}[\widehat{c}(\widetilde{e_{i}})]=0$ in the above frame
$\{\widetilde{e_{1}},\widetilde{e_{2}},\cdots,\widetilde{e_{n}}\}$. In terms of normal coordinates about $x_{0}$ one has:
$\sigma^{j}_{S(TM)}(x_{0})=0$ $e_{j}\big(c(e_{i})\big)(x_{0})=0$, $\Gamma^{k}(x_{0})=0$.
 \begin{lem}\cite{WW5} In the normal coordinates $U$ of $x_0$ in $M$,
 \begin{align}
 \partial_{x_j}(\sigma_{-2}(\widetilde{D}_{T}^{-2}))(x_0)=0,\partial_{x_j}(c(\xi))(x_0)=0.
\end{align}
\end{lem}
Substituting above results into the formula (2.18), we obtain
\begin{align}
 &\sigma_{-2m}(\widehat{c}(u)\widehat{c}(v) (d+\delta+\sqrt{-1}T_{1})^{1-2m})(x_{0})|_{|\xi|=1} \nonumber\\
 =&\widehat{c}(u)\widehat{c}(v) \Big\{|\xi|^{-2m}\sigma_0(d+\delta+T_{1})+\sum_{j=1}^{2m+2}\partial_{\xi_j}(|\xi|^{-2m})\partial_{x_j}(c(\xi))\nonumber\\
 &+
\Big[m\sigma_2(d+\delta+\sqrt{-1}T_{1})^{-m+1}\sigma_{-3}(d+\delta+\sqrt{-1}T_{1})^{2}
-\sqrt{-1} \sum_{k=0}^{m-2}\sum_{\mu=1}^{2m+2}\partial_{\xi_{\mu}}\sigma_{2}^{-m+k+1}((d+\delta+\sqrt{-1}T_{1})^2)\nonumber\\
&\times
\partial_{x_{\mu}}\sigma_{2}^{-1}((d+\delta+T_{1})^2)
(\sigma_2((d+\delta+\sqrt{-1}T_{1})^2))^{-k}\Big ]\sqrt{-1}c(\xi)\Big\}(x_{0})|_{|\xi|=1} \nonumber\\
=& \sum_{k<l }\sqrt{-1}T(\widetilde{e_k},\widetilde{e_l})\widehat{c}(u)\widehat{c}(v)
\widehat{c}(\widetilde{e_k})\widehat{c}(\widetilde{e_l})
 +m\sum_{i,j}g^{ij} \sum_{k<l }\sqrt{-1}T(\widetilde{e_k},\widetilde{e_l})\widehat{c}(u)\widehat{c}(v)\nonumber\\
 &\times\Big[c(\partial_{i})
\widehat{c}(\widetilde{e_k})\widehat{c}(\widetilde{e_l})
 + \widehat{c}(\widetilde{e_k})\widehat{c}(\widetilde{e_l}) c(\partial_{i})   \Big]\xi_{j}c(\xi).
\end{align}

Let $u=\sum_{i=1}^{n}u_{i}\widetilde{e_{i}}$, $v=\sum_{j=1}^{n}v_{j}\widetilde{e_{j}} $, $w=\sum_{l=1}^{n}w_{l}\widetilde{e_{l}} $,
where $\{\widetilde{e_{1}},\widetilde{e_{2}},\cdots,\widetilde{e_{n}}\}$ is the orthogonal basis about $g^{TM}$,
then  $c(u)= \sum_{i=1}^{n}u_{i}c(\widetilde{e_{i}})$,$c(v)= \sum_{j=1}^{n}v_{j}c(\widetilde{e_{j}})$,
$c(w)= \sum_{l=1}^{n}w_{l}c(\widetilde{e_{l}})$.
The following Lemmas of traces in terms of the Clifford action
is very efficient for solving the spectral form associated with the de-Rham Hodge type operator.
 \begin{lem}
The following identities hold:
 \begin{align}
 &{\rm{Tr}}\Big(\sum_{k<l }T(\widetilde{e_k},\widetilde{e_l})\widehat{c}(u)\widehat{c}(v)
\widehat{c}(\widetilde{e_k})\widehat{c}(\widetilde{e_l})\Big)(x_{0})|_{|\xi|=1}=-T(u,v){\rm{Tr}}(\rm{Id}).
\end{align}
\end{lem}
\begin{proof}
By the relation of the Clifford action and $ {\rm{Tr}}(AB)= {\rm{Tr}}(BA) $, we get
 \begin{align}
 &{\rm{Tr}}\Big(\sum_{k<l }T(\widetilde{e_k},\widetilde{e_l})\widehat{c}(u)\widehat{c}(v)
\widehat{c}(\widetilde{e_k})\widehat{c}(\widetilde{e_l})\Big)(x_{0}) \nonumber\\
 =&\sum^n_{k<l}\sum_{\alpha,\beta=1}^{n}T( \widetilde{e_k},\widetilde{e_l})u_{\alpha}v_{\beta}
{\rm{Tr}}\big( \widehat{c}(\widetilde{e_{\alpha}}) \widehat{c}(\widetilde{e_{\beta}}) \widehat{c}(\widetilde{e_{k}})
 \widehat{c}(\widetilde{e_{l}})\big)(x_{0})\nonumber\\
  =&\sum^n_{k<l}\sum_{\alpha,\beta=1}^{n}T( \widetilde{e_k},\widetilde{e_l})u_{\alpha}v_{\beta}
{\rm{Tr}}\Big( \widehat{c}(\widetilde{e_{\alpha}}) \big( -\widehat{c}(\widetilde{e_{k}})\widehat{c}(\widetilde{e_{\beta}})+2\delta_{k}^{\beta} \big)
 \widehat{c}(\widetilde{e_{l}})\Big)(x_{0})\nonumber\\
   =&-\sum^n_{k<l}\sum_{\alpha,\beta=1}^{n}T( \widetilde{e_k},\widetilde{e_l})u_{\alpha}v_{\beta}
{\rm{Tr}}\Big( \widehat{c}(\widetilde{e_{\alpha}})  \widehat{c}(\widetilde{e_{k}})\widehat{c}(\widetilde{e_{\beta}})
 \widehat{c}(\widetilde{e_{l}})\Big)(x_{0})\nonumber\\
   &+\sum^n_{k<l}\sum_{\alpha,\beta=1}^{n}T( \widetilde{e_k},\widetilde{e_l})u_{\alpha}v_{\beta}
2\delta_{k}^{\beta} {\rm{Tr}}\Big( \widehat{c}(\widetilde{e_{\alpha}})  \widehat{c}(\widetilde{e_{l}})\Big)(x_{0})\nonumber\\
   =&-\sum^n_{k<l}\sum_{\alpha,\beta=1}^{n}T( \widetilde{e_k},\widetilde{e_l})u_{\alpha}v_{\beta}
{\rm{Tr}}\Big( \big(-\widehat{c}(\widetilde{e_{k}})\widehat{c}(\widetilde{e_{\alpha}})+2\delta_{k}^{\alpha} \big)\widehat{c}(\widetilde{e_{\beta}})
 \widehat{c}(\widetilde{e_{l}})\Big)(x_{0})\nonumber\\
   &+\sum^n_{k<l}\sum_{\alpha,\beta=1}^{n}T( \widetilde{e_k},\widetilde{e_l})u_{\alpha}v_{\beta}
2\delta_{k}^{\beta} {\rm{Tr}}\Big( \widehat{c}(\widetilde{e_{\alpha}})  \widehat{c}(\widetilde{e_{l}})\Big)(x_{0})\nonumber\\
   =&\sum^n_{k<l}\sum_{\alpha,\beta=1}^{n}T( \widetilde{e_k},\widetilde{e_l})u_{\alpha}v_{\beta}
{\rm{Tr}}\Big(\widehat{c}(\widetilde{e_{k}})\widehat{c}(\widetilde{e_{\alpha}})\widehat{c}(\widetilde{e_{\beta}})
 \widehat{c}(\widetilde{e_{l}})\Big)(x_{0})\nonumber\\
 &-\sum^n_{k<l}\sum_{\alpha,\beta=1}^{n}T( \widetilde{e_k},\widetilde{e_l})u_{\alpha}v_{\beta}2\delta_{k}^{\alpha}
{\rm{Tr}}\Big(\widehat{c}(\widetilde{e_{\beta}})
 \widehat{c}(\widetilde{e_{l}})\Big)(x_{0})\nonumber\\
   &+\sum^n_{k<l}\sum_{\alpha,\beta=1}^{n}T( \widetilde{e_k},\widetilde{e_l})u_{\alpha}v_{\beta}
2\delta_{k}^{\beta} {\rm{Tr}}\Big( \widehat{c}(\widetilde{e_{\alpha}})  \widehat{c}(\widetilde{e_{l}})\Big)(x_{0})\nonumber\\
  =&-\sum^n_{k<l}\sum_{\alpha,\beta=1}^{n}T( \widetilde{e_k},\widetilde{e_l})u_{\alpha}v_{\beta}
{\rm{Tr}}\Big(\widehat{c}(\widetilde{e_{\alpha}})\widehat{c}(\widetilde{e_{\beta}})
\widehat{c}(\widetilde{e_{k}}) \widehat{c}(\widetilde{e_{l}})\Big)(x_{0})\nonumber\\
 &-\sum^n_{k<l}\sum_{\alpha,\beta=1}^{n}T( \widetilde{e_k},\widetilde{e_l})u_{\alpha}v_{\beta}2\delta_{k}^{\alpha}
{\rm{Tr}}\Big(\widehat{c}(\widetilde{e_{\beta}})
 \widehat{c}(\widetilde{e_{l}})\Big)(x_{0})\nonumber\\
   &+\sum^n_{k<l}\sum_{\alpha,\beta=1}^{n}T( \widetilde{e_k},\widetilde{e_l})u_{\alpha}v_{\beta}
2\delta_{k}^{\beta} {\rm{Tr}}\Big( \widehat{c}(\widetilde{e_{\alpha}})  \widehat{c}(\widetilde{e_{l}})\Big)(x_{0}).
\end{align}
By the generating relation, we see that the left hand side of (2.22) equals
 \begin{align}
&{\rm{Tr}}\Big(\sum_{k<l }T(\widetilde{e_k},\widetilde{e_l})\widehat{c}(u)\widehat{c}(v)
\widehat{c}(\widetilde{e_k})\widehat{c}(\widetilde{e_l})\Big)(x_{0})|_{|\xi|=1}\nonumber\\
=&-\sum^n_{k<l}\sum_{\alpha,\beta=1}^{n}T( \widetilde{e_k},\widetilde{e_l})u_{\alpha}v_{\beta} \delta_{k}^{\alpha}
{\rm{Tr}}\Big(\widehat{c}(\widetilde{e_{\beta}})
 \widehat{c}(\widetilde{e_{l}})\Big)(x_{0})\nonumber\\
   &+\sum^n_{k<l}\sum_{\alpha,\beta=1}^{n}T( \widetilde{e_k},\widetilde{e_l})u_{\alpha}v_{\beta}
 \delta_{k}^{\beta} {\rm{Tr}}\Big( \widehat{c}(\widetilde{e_{\alpha}})  \widehat{c}(\widetilde{e_{l}})\Big)(x_{0})\nonumber\\
=&\frac{1}{2!}\big(T(v,u)-T(u,v)\big){\rm{Tr}}(\rm{Id})\nonumber\\
=&-T(u,v){\rm{Tr}}(\rm{Id}).
\end{align}
Then the proof of the lemma is complete.
\end{proof}
We next make  argument under integral conditions to show  the calculation of traces in several special cases.
According to the integral formula of Case aII in \cite{YW1}, we have
 \begin{equation}
\int_{|\xi|=1}\xi_{i}\xi_{j}\sigma(\xi)=\frac{1}{n}\delta_{i}^{j}{\rm{vol}}(S^{n-1}).
\end{equation}

 \begin{lem}
In terms of the condition (2.24),
the following identities hold:
 \begin{align}
 &\int_{|\xi|=1}{\rm{Tr}}\Big(\sum_{i,j}g^{ij}\widehat{c}(u)\widehat{c}(v) c(\partial_{i})\sum_{k< l}
 T(\widetilde{e_k},\widetilde{e_l})
 \widehat{c}(\widetilde{e_k})\widehat{c}(\widetilde{e_l})\xi_{j}
 c(\xi) \Big)\sigma(\xi)(x_{0})|_{|\xi|=1}
 =T(u,v){\rm{Tr}}(\rm{Id}){\rm{vol}}(S^{2m-1});\\
 &\int_{|\xi|=1}{\rm{Tr}}\Big(\sum_{i,j}g^{ij}\widehat{c}(u)\widehat{c}(v)\sum_{k< l}
 T(\widetilde{e_k},\widetilde{e_l})
 \widehat{c}(\widetilde{e_k})\widehat{c}(\widetilde{e_l})c(\partial_{i})\xi_{j}
 c(\xi) \Big)\sigma(\xi)(x_{0})|_{|\xi|=1}
 =T(u,v){\rm{Tr}}(\rm{Id}){\rm{vol}}(S^{2m-1}).
\end{align}
\end{lem}
\begin{proof}
By the relation of the Clifford action and $ {\rm{Tr}}(AB)= {\rm{Tr}}(BA) $, in terms of the condition (2.24), we get
 \begin{align}
&\int_{|\xi|=1}{\rm{Tr}}\Big(\sum_{i,j}g^{ij}\widehat{c}(u)\widehat{c}(v) c(\partial_{i})\sum_{k< l}
 T(\widetilde{e_k},\widetilde{e_l}) \widehat{c}(\widetilde{e_k})\widehat{c}(\widetilde{e_l})\xi_{j}
 c(\xi) \Big)(x_{0})|_{|\xi|=1}\sigma(\xi)\nonumber\\
  =&\int_{|\xi|=1}\sum_{i }\sum_{k< l}T(\widetilde{e_k},\widetilde{e_l}){\rm{Tr}}\Big( \widehat{c}(u)\widehat{c}(v) c(\partial_{i})
 \widehat{c}(\widetilde{e_k})\widehat{c}(\widetilde{e_l})\xi_{i}
 \sum_{j }\xi_{j}c(\widetilde{e_j}) \Big)(x_{0})|_{|\xi|=1}\sigma(\xi)\nonumber\\
  =&\sum_{i }\sum_{k< l}T(\widetilde{e_k},\widetilde{e_l})\frac{1}{2m}\delta_{i}^{j}{\rm{vol}}(S^{2m-1})
  {\rm{Tr}}\Big(\widehat{c}(u)\widehat{c}(v)c(\widetilde{e_i}) \widehat{c}(\widetilde{e_k})
  \widehat{c}(\widetilde{e_l})c(\widetilde{e_j})\Big)(x_{0}) \nonumber\\
   =& \sum_{k< l}T(\widetilde{e_k},\widetilde{e_l}) {\rm{vol}}(S^{2m-1})
  {\rm{Tr}}\Big(c(\widetilde{e_i})\widehat{c}(u)\widehat{c}(v) \widehat{c}(\widetilde{e_k})
  \widehat{c}(\widetilde{e_l})c(\widetilde{e_i})\Big)(x_{0}) \nonumber\\
   =& -\sum_{k< l}T(\widetilde{e_k},\widetilde{e_l}){\rm{vol}}(S^{2m-1})
  {\rm{Tr}}\Big(\widehat{c}(u)\widehat{c}(v) \widehat{c}(\widetilde{e_k})
  \widehat{c}(\widetilde{e_l})\Big)(x_{0})=T(u,v){\rm{Tr}}(\rm{Id}){\rm{vol}}(S^{2m-1}).
\end{align}

On the other hand,
 \begin{align}
&\int_{|\xi|=1}{\rm{Tr}}\Big(\sum_{i,j}g^{ij}\widehat{c}(u)\widehat{c}(v)\sum_{k< l}
 T(\widetilde{e_k},\widetilde{e_l}) \widehat{c}(\widetilde{e_k})\widehat{c}(\widetilde{e_l}) c(\partial_{i})\xi_{j}
 c(\xi) \Big)(x_{0})|_{|\xi|=1}\sigma(\xi)\nonumber\\
  =&\int_{|\xi|=1}\sum_{i }\sum_{k< l}T(\widetilde{e_k},\widetilde{e_l}){\rm{Tr}}\Big( \widehat{c}(u)\widehat{c}(v)
 \widehat{c}(\widetilde{e_k})\widehat{c}(\widetilde{e_l})c(\partial_{i})\xi_{i}
 \sum_{j }\xi_{j}c(\widetilde{e_j}) \Big)(x_{0})|_{|\xi|=1}\sigma(\xi)\nonumber\\
  =&\sum_{i }\sum_{k< l}T(\widetilde{e_k},\widetilde{e_l})\frac{1}{2m}\delta_{i}^{j}{\rm{vol}}(S^{2m-1})
  {\rm{Tr}}\Big(\widehat{c}(u)\widehat{c}(v)\widehat{c}(\widetilde{e_k})
  \widehat{c}(\widetilde{e_l})c(\widetilde{e_i}) c(\widetilde{e_j})\Big)(x_{0}) \nonumber\\
   =& -\sum_{k< l}T(\widetilde{e_k},\widetilde{e_l}){\rm{vol}}(S^{2m-1})
  {\rm{Tr}}\Big(\widehat{c}(u)\widehat{c}(v) \widehat{c}(\widetilde{e_k})
  \widehat{c}(\widetilde{e_l})\Big)(x_{0})=T(u,v){\rm{Tr}}(\rm{Id}){\rm{vol}}(S^{2m-1}).
\end{align}
Then the proof of the lemma is complete.
\end{proof}
From (2.14), Lemma 2.4 and Lemma 2.5, then  the proof of Theorem 1.1 is complete.

\section{Trilinear functional  for the  de-Rham Hodge type operator with 3-form torsion}

Let $M$ be a $n$-dimensional ($n\geq 3$) oriented compact Riemannian manifold with a Riemannian metric $g^{M}$.
 And let $\nabla^L$ be the Levi-Civita connection about $g^M$. In the local coordinates $\{x_i; 1\leq i\leq n\}$ and the
fixed orthonormal frame $\{\widetilde{e_1},\cdots,\widetilde{e_n}\}$, the connection matrix $(\omega_{s,t})$ is defined by
\begin{equation}
\nabla^L(\widetilde{e_1},\cdots,\widetilde{e_n})= (\widetilde{e_1},\cdots,\widetilde{e_n})(\omega_{s,t}).
\end{equation}
The Levi-Civita connection
$\nabla^L: \Gamma(TM)\rightarrow \Gamma(T^{*}M\otimes TM)$ on $M$ induces a connection
$\nabla^{S}: \Gamma(S)\rightarrow \Gamma(T^{*}M\otimes S).$ By adding a additional torsion term $T\in\Omega^{1}(M,{\rm{ End}}(TM))$ we
obtain a new covariant derivative
 \begin{equation}
\langle \widetilde{\nabla}^{T}_{\widetilde{e_i}}\widetilde{e_s}, \widetilde{e_t} \rangle=
\langle \nabla^{L}_{\widetilde{e_i}}\widetilde{e_s}, \widetilde{e_t} \rangle+T(\widetilde{e_i},\widetilde{e_s},\widetilde{e_t}),
\end{equation}
where the 3-form $T=\sum _{1\leq\alpha<\beta<\gamma\leq n}
 T_{\alpha \beta \gamma}c(e_{\alpha})c(e_{\beta})c(e_{\gamma})$ denotes the `lifted' torsion term.
 The  metric connection $\widetilde{\nabla} $ is in fact compatible with the Riemannian metric $g$ and therefore also induces a
 connection on $ \wedge^{*}(T^{*}M)$\cite{Y, GHV}, and the de-Rham Hodge type operator with torsion as follows.
\begin{align}
(d+\delta)_{T_{2}}=&\sum^n_{i=1}c(\widetilde{e_i})\bigg[\widetilde{e_i}-\frac{1}{4}
\sum_{s,t}\langle \nabla^{T_{2}}_{\widetilde{e_i}}\widetilde{e_s}, \widetilde{e_t} \rangle
\Big[\widehat{c}(\widetilde{e_s})\widehat{c}(\widetilde{e_t})-c(\widetilde{e_s})c(\widetilde{e_t})\Big]\bigg]\nonumber\\
=&\sum^n_{i=1}c(\widetilde{e_i})\bigg[\widetilde{e_i}-\frac{1}{4}\sum_{s,t}
\Big(\langle \nabla^{L}_{\widetilde{e_i}}\widetilde{e_s}, \widetilde{e_t} \rangle+T(\widetilde{e_i},\widetilde{e_s},\widetilde{e_t}) \Big)
\Big[\widehat{c}(\widetilde{e_s})\widehat{c}(\widetilde{e_t})-c(\widetilde{e_s})c(\widetilde{e_t})\Big]\bigg]\nonumber\\
=& \sum^n_{i=1}c(\widetilde{e_i})\bigg[\widetilde{e_i}-\frac{1}{4}\sum_{s,t}\omega_{s,t}
(\widetilde{e_i})\Big(\widehat{c}(\widetilde{e_s})\widehat{c}(\widetilde{e_t})-c(\widetilde{e_s})c(\widetilde{e_t})\Big)\bigg]
-\frac{1}{4}\sum_{i,s,t}T(\widetilde{e_i},\widetilde{e_s},\widetilde{e_t})
c(\widetilde{e_i})\widehat{c}(\widetilde{e_s})\widehat{c}(\widetilde{e_t})\nonumber\\
&+\frac{1}{4}\sum_{i,s,t}T(\widetilde{e_i},\widetilde{e_s},\widetilde{e_t})
c(\widetilde{e_i})c(\widetilde{e_s})c(\widetilde{e_t})\nonumber\\
=& \sum^n_{i=1}c(\widetilde{e_i})\bigg[\widetilde{e_i}+\frac{1}{4}\sum_{s,t}\omega_{s,t}
(\widetilde{e_i})[\widehat{c}(\widetilde{e_s})\widehat{c}(\widetilde{e_t})-c(\widetilde{e_s})c(\widetilde{e_t})]\bigg]
+\frac{3}{2}\sum_{i<s<t}T(\widetilde{e_i},\widetilde{e_s},\widetilde{e_t})
c(\widetilde{e_i})c(\widetilde{e_s})c(\widetilde{e_t})\nonumber\\
&-\frac{1}{4}\sum^n_{i\neq s\neq t}T(\widetilde{e_i},\widetilde{e_s},\widetilde{e_t})
c(\widetilde{e_i})\widehat{c}(\widetilde{e_s})\widehat{c}(\widetilde{e_t}).
\end{align}
Using an explicit formula for $d+\delta$, we can
reformulate $(d+\delta)_{T_{2}}$ as follows.
\begin{defn}
The de-Rham Hodge type operator  is the first order differential operator on $ \Gamma(M,\wedge^{*}(T^{*}M))$  given by the formula
 \begin{align}
(d+\delta)_{T_{2}}=&d+\delta
+\frac{3}{2}\sum_{i<s<t}T(\widetilde{e_i},\widetilde{e_s},\widetilde{e_t})
c(\widetilde{e_i})c(\widetilde{e_s})c(\widetilde{e_t})
-\frac{1}{4}\sum^n_{i\neq s\neq t}T(\widetilde{e_i},\widetilde{e_s},\widetilde{e_t})
c(\widetilde{e_i})\widehat{c}(\widetilde{e_s})\widehat{c}(\widetilde{e_t}).
\end{align}
\end{defn}
In this section we want to consider the trilinear functional  for the  de-Rham Hodge type operator with torsion.
To avoid technique terminology we only state our results for
compact oriented Riemannian manifold of even dimension $n=2m$ by using the trace of  de-Rham Hodge type operator and the noncommutative residue density.
\begin{defn}\cite{DSZ}
For $(d+\delta)_{T_{2}}$ given by (3.4), the trilinear Clifford multiplication by functional of differential one-forms
$c(u),c(v),c(w)$
\begin{align}
\mathscr{T}_{2}(c(u),c(v),c(w))={\rm Wres}\big(c(u)c(v)c(w)(d+\delta)_{T_{2}}^{-n+1}\big)
\end{align}
is called torsion functional $\mathscr{T}_{2}$ associated with the de-Rham Hodge type operator.
\end{defn}
\begin{defn}\cite{DSZ}
For $(d+\delta)_{T_{2}}$ given by (3.4), the trilinear Clifford multiplication by functional of differential one-forms
$c(u),\widehat{c}(v),\widehat{c}(w)$
\begin{align}
\mathscr{T}_{3}(c(u),\widehat{c}(v),\widehat{c}(w))={\rm Wres}\big(c(u)\widehat{c}(v)\widehat{c}(w)(d+\delta)_{T_{2}}^{-n+1}\big)
\end{align}
is called torsion functional $\mathscr{T}_{3}$ associated with the de-Rham Hodge type operator.
\end{defn}
The following noncommutative residue formula \cite{DSZ} plays a key role in our proof
of the spectral torsion  associated with the de-Rham Hodge type operator.
\begin{align}
{\rm Wres}\big(c(u)c(v)c(w)(d+\delta)_{T_{2}}^{-2m+1}\big)=\int_M\int_{|\xi|=1}{\rm Tr}_{\wedge^{*}(T^{*}M)}
[\sigma_{-2m}\big(c(u)c(v)c(w)(d+\delta)_{T_{2}}^{-2m+1}\big)]\sigma(\xi){\rm d}x.
\end{align}

\subsection{The proof of Theorem 1.2}
For any fixed point $x_0\in M$, we can choose the normal coordinates
$U$ of $x_0$ in $M$ and compute $\sigma_{-n}\big(c(u)c(v)c(w)(d+\delta)_{T_{2}}^{-n+1}\big)$.
Let $c(\widetilde{e_j})=\epsilon (\widetilde{e_j*})-\iota (\widetilde{e_j*});\
\widehat{c}(\widetilde{e_j})=\epsilon (\widetilde{e_j*})+\iota (\widetilde{e_j*})$ act on $\Gamma(M,\wedge^{*}(T^{*}M))$.
Then we have $ \partial_{x_j}[c(\widetilde{e_{i}})]=\partial_{x_j}[\widehat{c}(\widetilde{e_{i}})]=0$ in the above frame
$\{\widetilde{e_{1}},\widetilde{e_{2}},\cdots,\widetilde{e_{n}}\}$. In terms of normal coordinates about $x_{0}$ one has:
$\sigma^{j}_{S(TM)}(x_{0})=0$ $e_{j}\big(c(e_{i})\big)(x_{0})=0$, $\Gamma^{k}(x_{0})=0$.
Substituting above results into the formula (2.18), we obtain
\begin{align}
 &\sigma_{-2m}(c(u)c(v)c(w)(d+\delta)_{T_{2}}^{1-2m})(x_{0})|_{|\xi|=1} \nonumber\\
 =&c(u)c(v)c(w)\Big\{|\xi|^{-2m}\sigma_0((d+\delta)_{T_{2}})+\sum_{j=1}^{2m+2}\partial_{\xi_j}(|\xi|^{-2m})\partial_{x_j}c(\xi)+
\Big[m\sigma_2((d+\delta)_{T_{2}})^{-m+1}\sigma_{-3}((d+\delta)_{T_{2}}^{-2})  \nonumber\\
&-\sqrt{-1} \sum_{k=0}^{m-2}\sum_{\mu=1}^{2m+2}
\partial_{\xi_{\mu}}\sigma_{2}^{-m+k+1}((d+\delta)_{T_{2}}^2)
\partial_{x_{\mu}}\sigma_{2}^{-1}((d+\delta)_{T_{2}}^2)(\sigma_2((d+\delta)_{T_{2}}^2))^{-k}\Big ]\sqrt{-1}c(\xi)\Big\}(x_{0})|_{|\xi|=1} \nonumber\\
=&c(u)c(v)c(w)\Big(\frac{3}{2}\sum_{k<s<t}T(\widetilde{e_k},\widetilde{e_s},\widetilde{e_t})
c(\widetilde{e_k})c(\widetilde{e_s})c(\widetilde{e_t})
-\frac{1}{4}\sum^n_{k\neq s\neq t}T(\widetilde{e_k},\widetilde{e_s},\widetilde{e_t})
c(\widetilde{e_k})\widehat{c}(\widetilde{e_s})\widehat{c}(\widetilde{e_t})\Big)\nonumber\\
&+m\sum_{i,j}g^{ij}c(u)c(v)c(w)\Big[c(\partial_{i})\Big(\frac{3}{2}\sum_{k<s<t}T(\widetilde{e_k},\widetilde{e_s},\widetilde{e_t})
c(\widetilde{e_k})c(\widetilde{e_s})c(\widetilde{e_t})
-\frac{1}{4}\sum^n_{k\neq s\neq t}T(\widetilde{e_k},\widetilde{e_s},\widetilde{e_t})
c(\widetilde{e_k})\widehat{c}(\widetilde{e_s})\widehat{c}(\widetilde{e_t})\Big) \nonumber\\
&+\Big(\frac{3}{2}\sum_{k<s<t}T(\widetilde{e_k},\widetilde{e_s},\widetilde{e_t})
c(\widetilde{e_k})c(\widetilde{e_s})c(\widetilde{e_t})
-\frac{1}{4}\sum^n_{k\neq s\neq t}T(\widetilde{e_k},\widetilde{e_s},\widetilde{e_t})
c(\widetilde{e_k})\widehat{c}(\widetilde{e_s})\widehat{c}(\widetilde{e_t})\Big)c(\partial_{i})   \Big]\xi_{j}c(\xi).
\end{align}

Let $u=\sum_{i=1}^{n}u_{i}\widetilde{e_{i}}$, $v=\sum_{j=1}^{n}v_{j}\widetilde{e_{j}} $, $w=\sum_{l=1}^{n}w_{l}\widetilde{e_{l}} $,
where $\{\widetilde{e_{1}},\widetilde{e_{2}},\cdots,\widetilde{e_{n}}\}$ is the orthogonal basis about $g^{TM}$,
then  $c(u)= \sum_{i=1}^{n}u_{i}c(\widetilde{e_{i}})$,$c(v)= \sum_{j=1}^{n}v_{j}c(\widetilde{e_{j}})$,
$c(w)= \sum_{l=1}^{n}w_{l}c(\widetilde{e_{l}})$.
The following Lemmas of traces in terms of the Clifford action
is very efficient for solving the spectral torsion  associated with the de-Rham Hodge type operator.
 \begin{lem}
The following identities hold:
 \begin{align}
 &{\rm{Tr}}\Big(c(u)c(v)c(w)\sum_{i<s<t}T(\widetilde{e_i},\widetilde{e_s},\widetilde{e_t})
c(\widetilde{e_i})c(\widetilde{e_s})c(\widetilde{e_t})\Big)=T(u,v,w){\rm{Tr}}(\rm{Id});\\
&{\rm{Tr}}\Big(c(u)c(v)c(w)\sum^n_{i\neq s\neq t}T(\widetilde{e_i},\widetilde{e_s},\widetilde{e_t})
c(\widetilde{e_i})\widehat{c}(\widetilde{e_s})\widehat{c}(\widetilde{e_t})\Big)=0.
\end{align}
\end{lem}
\begin{proof}
By the relation of the Clifford action and $ {\rm{Tr}}(AB)= {\rm{Tr}}(BA) $, we get
 \begin{align}
 &{\rm{Tr}}\Big(c(u)c(v)c(w)\sum^n_{j< k< l}T(\widetilde{e_j},\widetilde{e_k},\widetilde{e_l})
c(\widetilde{e_j})c(\widetilde{e_k})c(\widetilde{e_l})\Big)\nonumber\\
 =&\sum^n_{j< k< l}\sum_{\alpha,\beta,\gamma=1}^{n}T(\widetilde{e_j},\widetilde{e_k},\widetilde{e_l})u_{\alpha}v_{\beta}w_{\gamma}
{\rm{Tr}}\big(c(\widetilde{e_{\alpha}})c(\widetilde{e_{\beta}})c(\widetilde{e_{\gamma}})c(\widetilde{e_{j}})c(\widetilde{e_{k}})
c(\widetilde{e_{l}})\big)\nonumber\\
   =&\sum^n_{j< k< l}\sum_{\alpha,\beta,\gamma=1}^{n}T(\widetilde{e_j},\widetilde{e_k},\widetilde{e_l})u_{\alpha}v_{\beta}w_{\gamma}
{\rm{Tr}}\big(c(\widetilde{e_{\alpha}})c(\widetilde{e_{\beta}})c(\widetilde{e_{\gamma}})
c(\widetilde{e_{l}})c(\widetilde{e_{j}})c(\widetilde{e_{k}})\big)\nonumber\\
  =&\sum^n_{j< k< l}\sum_{\alpha,\beta,\gamma=1}^{n}T(\widetilde{e_j},\widetilde{e_k},\widetilde{e_l})u_{\alpha}v_{\beta}w_{\gamma}
{\rm{Tr}}\Big(c(\widetilde{e_{\alpha}})c(\widetilde{e_{\beta}})\big(-c(\widetilde{e_{l}})c(\widetilde{e_{\gamma}})
-2\delta_{\gamma}^{l}\big)c(\widetilde{e_{j}})c(\widetilde{e_{k}})\Big)\nonumber\\
=&\cdots   \nonumber\\
=&-\sum^n_{j< k< l}\sum_{\alpha,\beta,\gamma=1}^{n}T(\widetilde{e_j},\widetilde{e_k},\widetilde{e_l})u_{\alpha}v_{\beta}w_{\gamma}
{\rm{Tr}}\big(c(\widetilde{e_{\alpha}})c(\widetilde{e_{\beta}})c(\widetilde{e_{\gamma}})c(\widetilde{e_{j}})
c(\widetilde{e_{k}})c(\widetilde{e_{l}})\big)\nonumber\\
&+2\sum_{1 \leq  j,k,l\leq n}\sum_{\alpha,\beta,\gamma=1}^{n}T(\widetilde{e_j},\widetilde{e_k},\widetilde{e_l})u_{\alpha}v_{\beta}w_{\gamma}
\Big(\delta_{\alpha}^{l}\delta_{\beta}^{j }\delta_{\gamma}^{k}-\delta_{\alpha}^{l}\delta_{\beta}^{k}\delta_{\gamma}^{j}
+\delta_{\alpha}^{l}\delta_{\beta}^{\gamma}\delta_{j}^{k}\nonumber\\
&+\delta_{\gamma}^{l}\delta_{\alpha}^{j}\delta_{\beta}^{k}-\delta_{\gamma}^{l}\delta_{\alpha}^{k}\delta_{\beta}^{j}
+\delta_{\gamma}^{l}\delta_{\alpha}^{\beta}\delta_{j}^{k}
-\delta_{\beta}^{l}\delta_{\alpha}^{j}\delta_{\gamma}^{k}+\delta_{\beta}^{l}\delta_{\alpha}^{k}\delta_{\gamma}^{j}
+\delta_{\beta}^{l}\delta_{\alpha}^{\gamma}\delta_{j}^{k} \Big){\rm{Tr}}(\rm{Id}).
\end{align}
By the generating relation, we see that the left hand side of (3.11) equals
 \begin{align}
 {\rm{Tr}}\Big(c(u)c(v)c(w)\sum^n_{i< s< t}T(\widetilde{e_i},\widetilde{e_s},\widetilde{e_t})
c(\widetilde{e_i})c(\widetilde{e_s})c(\widetilde{e_t})\Big)
=T(u,v,w){\rm{Tr}}(\rm{Id}).
\end{align}
In the same way we get
 \begin{align}
 &{\rm{Tr}}\Big(c(u)c(v)c(w)\sum^n_{i\neq s\neq t}T(\widetilde{e_i},\widetilde{e_s},\widetilde{e_t})
c(\widetilde{e_i})\widehat{c}(\widetilde{e_s})\widehat{c}(\widetilde{e_t})\Big) \nonumber\\
=&-{\rm{Tr}}\Big(\sum^n_{i\neq s\neq t}T(\widetilde{e_i},\widetilde{e_s},\widetilde{e_t})c(u)c(v)c(w)
\widehat{c}(\widetilde{e_s})c(\widetilde{e_i})\widehat{c}(\widetilde{e_t})\Big)\nonumber\\
=&\cdots \nonumber\\
=&{\rm{Tr}}\Big(\sum^n_{i\neq s\neq t}T(\widetilde{e_i},\widetilde{e_s},\widetilde{e_t})c(u)c(v)c(w)
c(\widetilde{e_i})\widehat{c}(\widetilde{e_t})\widehat{c}(\widetilde{e_s})\Big)\nonumber\\
=&{\rm{Tr}}\Big(\sum^n_{i\neq s\neq t}T(\widetilde{e_i},\widetilde{e_s},\widetilde{e_t})c(u)c(v)c(w)
c(\widetilde{e_i})\big(-\widehat{c}(\widetilde{e_s})\widehat{c}(\widetilde{e_t})+2\delta_{s}^{t} \big)\Big)\nonumber\\
=&-{\rm{Tr}}\Big(\sum^n_{i\neq s\neq t}T(\widetilde{e_i},\widetilde{e_s},\widetilde{e_t})c(u)c(v)c(w)
c(\widetilde{e_i})\widehat{c}(\widetilde{e_s})\widehat{c}(\widetilde{e_t})\Big).
\end{align}
Then the proof of the lemma is complete.
\end{proof}
We next make  argument under integral conditions to show  the calculation of traces in several special cases.
 \begin{lem}
In terms of the condition (2.24),
the following identities hold:
 \begin{align}
 &\int_{|\xi|=1}{\rm{Tr}}\Big(\sum_{i,j}g^{ij}c(u)c(v)c(w)c(\partial_{i})\sum_{l\neq s\neq t}T(\widetilde{e_l},\widetilde{e_s},\widetilde{e_t})
 c(\widetilde{e_l})\widehat{c}(\widetilde{e_s})\widehat{c}(\widetilde{e_t})\xi_{j}
 c(\xi) \Big)\sigma(\xi)
 =0;\\
 &\int_{|\xi|=1}{\rm{Tr}}\Big(\sum_{i,j}g^{ij}c(u)c(v)c(w)\sum_{l\neq s\neq t}T(\widetilde{e_l},\widetilde{e_s},\widetilde{e_t})
 c(\widetilde{e_l})\widehat{c}(\widetilde{e_s})\widehat{c}(\widetilde{e_t})c(\partial_{i})\xi_{j}c(\xi) \Big)\sigma(\xi)
 =0.
\end{align}
\end{lem}
\begin{proof}
By the relation of the Clifford action and $ {\rm{Tr}}(AB)= {\rm{Tr}}(BA) $, in terms of the condition (2.24), we get
 \begin{align}
&{\rm{Tr}}\Big(\sum_{l\neq s\neq t}c(u)c(v)c(w)c(\widetilde{e_i})c(\widetilde{e_l})\widehat{c}(\widetilde{e_s})\widehat{c}(\widetilde{e_t})
c(\widetilde{e_i}) \Big)\nonumber\\
  =&\sum_{h,j,k}\sum_{l\neq s\neq t}{\rm{Tr}}\Big(c(\widetilde{e_h})c(\widetilde{e_j})c(\widetilde{e_k})
  c(\widetilde{e_i})c(\widetilde{e_l})\widehat{c}(\widetilde{e_s})\widehat{c}(\widetilde{e_t})
c(\widetilde{e_i}) \Big)\nonumber\\
 =&\sum_{h,j,k}\sum_{l\neq s\neq t}{\rm{Tr}}\Big(c(\widetilde{e_h})c(\widetilde{e_j})
 \big(-c(\widetilde{e_i})c(\widetilde{e_k})-2\delta_{k}^{i}\big)
  c(\widetilde{e_l})\widehat{c}(\widetilde{e_s})\widehat{c}(\widetilde{e_t})
c(\widetilde{e_i}) \Big)\nonumber\\
=&-\sum_{h,j,k}\sum_{l\neq s\neq t}{\rm{Tr}}\Big(c(\widetilde{e_h})c(\widetilde{e_j})
c(\widetilde{e_i})c(\widetilde{e_k})
  c(\widetilde{e_l})\widehat{c}(\widetilde{e_s})\widehat{c}(\widetilde{e_t})
c(\widetilde{e_i}) \Big)\nonumber\\
&-2\sum_{h,j,k}\sum_{l\neq s\neq t}\delta_{k}^{i}{\rm{Tr}}\Big(c(\widetilde{e_h})c(\widetilde{e_j})
  c(\widetilde{e_l})\widehat{c}(\widetilde{e_s})\widehat{c}(\widetilde{e_t})
c(\widetilde{e_i}) \Big)\nonumber\\
=&-\sum_{h,j,k}\sum_{l\neq s\neq t}{\rm{Tr}}\Big(c(\widetilde{e_h})
\big(-c(\widetilde{e_i})c(\widetilde{e_j})-2\delta_{j}^{i}\big)
c(\widetilde{e_k})
  c(\widetilde{e_l})\widehat{c}(\widetilde{e_s})\widehat{c}(\widetilde{e_t})
c(\widetilde{e_i}) \Big)\nonumber\\
&-2\sum_{h,j,k}\sum_{l\neq s\neq t}\delta_{k}^{i}{\rm{Tr}}\Big(c(\widetilde{e_h})c(\widetilde{e_j})
  c(\widetilde{e_l})\widehat{c}(\widetilde{e_s})\widehat{c}(\widetilde{e_t})
c(\widetilde{e_i}) \Big)\nonumber\\
=&\cdots\nonumber\\
=&-\sum_{h,j,k}\sum_{l\neq s\neq t}{\rm{Tr}}\Big(c(\widetilde{e_i})c(\widetilde{e_h})
c(\widetilde{e_j})
c(\widetilde{e_k}) c(\widetilde{e_l})\widehat{c}(\widetilde{e_s})\widehat{c}(\widetilde{e_t})
c(\widetilde{e_i}) \Big)\nonumber\\
&-2\sum_{h,j,k}\sum_{l\neq s\neq t}\delta_{h}^{i}{\rm{Tr}}\Big(
c(\widetilde{e_j})
  c(\widetilde{e_l})\widehat{c}(\widetilde{e_s})\widehat{c}(\widetilde{e_t})
c(\widetilde{e_i}) \Big)\nonumber\\
&-2\sum_{h,j,k}\sum_{l\neq s\neq t}\delta_{j}^{i}{\rm{Tr}}\Big(c(\widetilde{e_h})
  c(\widetilde{e_l})\widehat{c}(\widetilde{e_s})\widehat{c}(\widetilde{e_t})
c(\widetilde{e_i}) \Big)\nonumber\\
&-2\sum_{h,j,k}\sum_{l\neq s\neq t}\delta_{k}^{i}{\rm{Tr}}\Big(c(\widetilde{e_h})c(\widetilde{e_j})
  c(\widetilde{e_l})\widehat{c}(\widetilde{e_s})\widehat{c}(\widetilde{e_t})
c(\widetilde{e_i}) \Big).
\end{align}
Then we obtain
\begin{align}
 &\int_{|\xi|=1}{\rm{Tr}}\Big(\sum_{i,j}g^{ij}c(u)c(v)c(w)c(\partial_{i})\sum_{l\neq s\neq t}T(\widetilde{e_l},\widetilde{e_s},\widetilde{e_t})
 c(\widetilde{e_l})\widehat{c}(\widetilde{e_s})\widehat{c}(\widetilde{e_t})\xi_{j}c(\xi) \Big)\sigma(\xi)\nonumber\\
 =&\int_{|\xi|=1}{\rm{Tr}}\Big(\sum_{i,j}g^{ij}c(u)c(v)c(w) c(\partial_{i})\sum_{l\neq s\neq t}T(\widetilde{e_l},\widetilde{e_s},\widetilde{e_t})
 c(\widetilde{e_l})\widehat{c}(\widetilde{e_s})\widehat{c}(\widetilde{e_t})
\xi_{j}\sum_{h}\xi_{h}c(e_{h}) \Big)\sigma(\xi)\nonumber\\
   =&\int_{|\xi|=1} {\rm{Tr}}\Big(\sum_{i,h}\xi_{i}\xi_{h}  c(u)c(v)c(w)c(\partial_{i})\sum_{l\neq s\neq t}T(\widetilde{e_l},\widetilde{e_s},\widetilde{e_t})
    c(\widetilde{e_l})\widehat{c}(\widetilde{e_s})\widehat{c}(\widetilde{e_t})
c(e_{h}) \Big)\sigma(\xi)\nonumber\\
  =&\sum_{i}\delta _{i}^{h}\frac{1}{n}{\rm{vol}}(S^{n-1})\sum_{l\neq s\neq t}T(\widetilde{e_l},\widetilde{e_s},\widetilde{e_t})
  {\rm{Tr}}\Big(c(u)c(v)c(w)c(\partial_{i})
   c(\widetilde{e_l})\widehat{c}(\widetilde{e_s})\widehat{c}(\widetilde{e_t})
c(e_{h}) \Big)\nonumber\\
  =&0.
\end{align}

On the other hand,
 \begin{align}
 &\int_{|\xi|=1}{\rm{Tr}}\Big(\sum_{i,j}g^{ij}c(u)c(v)c(w)\sum_{l\neq s\neq t}T(\widetilde{e_l},\widetilde{e_s},\widetilde{e_t})
 c(\widetilde{e_l})\widehat{c}(\widetilde{e_s})\widehat{c}(\widetilde{e_t})c(\partial_{i})\xi_{j}c(\xi) \Big)\sigma(\xi)\nonumber\\
 =&\int_{|\xi|=1}{\rm{Tr}}\Big(\sum_{i,j}g^{ij}c(u)c(v)c(w)\sum_{l\neq s\neq t}T(\widetilde{e_l},\widetilde{e_s},\widetilde{e_t})
 c(\widetilde{e_l})\widehat{c}(\widetilde{e_s})\widehat{c}(\widetilde{e_t})
 c(\partial_{i})\xi_{j}\sum_{h}\xi_{h}c(e_{h}) \Big)\sigma(\xi)\nonumber\\
   =&\int_{|\xi|=1} {\rm{Tr}}\Big(\sum_{i,h}\xi_{i}\xi_{h}  c(u)c(v)c(w)\sum_{l\neq s\neq t}T(\widetilde{e_l},\widetilde{e_s},\widetilde{e_t})
    c(\widetilde{e_l})\widehat{c}(\widetilde{e_s})\widehat{c}(\widetilde{e_t})
  c(\partial_{i})c(e_{h}) \Big)\sigma(\xi)\nonumber\\
  =&\sum_{i}\delta _{i}^{h}\frac{1}{n}{\rm{vol}}(S^{n-1})
  {\rm{Tr}}\Big(c(u)c(v)c(w)\sum_{l\neq s\neq t}T(\widetilde{e_l},\widetilde{e_s},\widetilde{e_t})
   c(\widetilde{e_l})\widehat{c}(\widetilde{e_s})\widehat{c}(\widetilde{e_t})
  c(\partial_{i})c(e_{l}) \Big)\nonumber\\
   =&-{\rm{vol}}(S^{n-1}){\rm{Tr}}\Big( c(u)c(v)c(w)\sum_{l\neq s\neq t}T(\widetilde{e_l},\widetilde{e_s},\widetilde{e_t})
    c(\widetilde{e_l})\widehat{c}(\widetilde{e_s})\widehat{c}(\widetilde{e_t}) \Big)\nonumber\\
 =&0.
\end{align}
Then the proof of the lemma is complete.
\end{proof}

 \begin{lem} \cite{WW5}
The following identities hold:
 \begin{align}
 &\int_{|\xi|=1}{\rm{Tr}}\Big(\sum_{i,j}g^{ij}c(u)c(v)c(w)
 \sum_{l<s<t}T(\widetilde{e_l},\widetilde{e_s},\widetilde{e_t})c(\widetilde{e_l})c(\widetilde{e_s})c(\widetilde{e_t})
 c(\partial_{i})\xi_{j}c(\xi) \Big)\sigma(\xi)\nonumber\\
 =&-T(u,v,w){\rm{Tr}}(\rm{Id}){\rm{vol}}(S^{n-1});\\
 &\int_{|\xi|=1}{\rm{Tr}}\Big(\sum_{i,j}g^{ij}c(u)c(v)c(w)c(\partial_{i})
 \sum_{l<s<t}T(\widetilde{e_l},\widetilde{e_s},\widetilde{e_t})c(\widetilde{e_l})c(\widetilde{e_s})c(\widetilde{e_t})
 \xi_{j}c(\xi) \Big)
 \sigma(\xi)\nonumber\\
 =&-5T(u,v,w){\rm{Tr}}(\rm{Id}){\rm{vol}}(S^{n-1}).
\end{align}
\end{lem}
It follows from Lemma 3.4, Lemma 3.5 and Lemma 3.6, we obtain
\begin{align}
&\int_{|\xi|=1}{\rm Tr}_{\wedge^{*}(T^{*}M)}
[\sigma_{-2m}\big(c(u)c(v)c(w)(d+\delta)_{T_{2} }^{-2m+1}\big)](x_{0})|_{|\xi|=1}\sigma(\xi) \nonumber\\
=&\int_{|\xi|=1}{\rm Tr}_{\wedge^{*}(T^{*}M)}\Big[
c(u)c(v)c(w)\Big(\frac{3}{2}\sum_{i<s<t}T(\widetilde{e_i},\widetilde{e_s},\widetilde{e_t})
c(\widetilde{e_i})c(\widetilde{e_s})c(\widetilde{e_t})
 -\frac{1}{4}\sum^n_{i\neq s\neq t}T(\widetilde{e_i},\widetilde{e_s},\widetilde{e_t})
c(\widetilde{e_i})\widehat{c}(\widetilde{e_s})\widehat{c}(\widetilde{e_t})\Big)\nonumber\\
&+m\sum_{i,j}g^{ij}c(u)\widehat{c}(v)\widehat{c}(w)\Big[c(\partial_{i})\Big(\frac{3}{2}\sum_{i<s<t}T(\widetilde{e_i},\widetilde{e_s},\widetilde{e_t})
c(\widetilde{e_i})c(\widetilde{e_s})c(\widetilde{e_t})
 -\frac{1}{4}\sum^n_{i\neq s\neq t}T(\widetilde{e_i},\widetilde{e_s},\widetilde{e_t})
c(\widetilde{e_i})\widehat{c}(\widetilde{e_s})\widehat{c}(\widetilde{e_t})\Big) \nonumber\\
&+\Big(\frac{3}{2}\sum_{i<s<t}T(\widetilde{e_i},\widetilde{e_s},\widetilde{e_t})
c(\widetilde{e_i})c(\widetilde{e_s})c(\widetilde{e_t})
 -\frac{1}{4}\sum^n_{i\neq s\neq t}T(\widetilde{e_i},\widetilde{e_s},\widetilde{e_t})
c(\widetilde{e_i})\widehat{c}(\widetilde{e_s})\widehat{c}(\widetilde{e_t})\Big)c(\partial_{i})   \Big]\xi_{j}c(\xi)\Big](x_{0})|_{|\xi|=1}\sigma(\xi)
 \nonumber\\
 =&2^{2m-1}(3-18m)T(u,v,w)
{\rm{vol}}(S^{2m-1}).
\end{align}
Substituting (3.21) into (3.7) leads to the desired equality, and the proof of Theorem 1.2 is complete.

\subsection{The proof of Theorem 1.3}

Also, from (2.18) we have
\begin{align}
 &\sigma_{-2m}(c(u)\widehat{c}(v)\widehat{c}(w)(d+\delta)_{T_{2}}^{1-2m})(x_{0})|_{|\xi|=1} \nonumber\\
=&c(u)\widehat{c}(v)\widehat{c}(w)\Big(\frac{3}{2}\sum_{i<s<t}T(\widetilde{e_i},\widetilde{e_s},\widetilde{e_t})
c(\widetilde{e_i})c(\widetilde{e_s})c(\widetilde{e_t})
 -\frac{1}{4}\sum^n_{i\neq s\neq t}T(\widetilde{e_i},\widetilde{e_s},\widetilde{e_t})
c(\widetilde{e_i})\widehat{c}(\widetilde{e_s})\widehat{c}(\widetilde{e_t}))\Big)\nonumber\\
&+m\sum_{i,j}g^{ij}c(u)\widehat{c}(v)\widehat{c}(w)\Big[c(\partial_{i})\Big(\frac{3}{2}\sum_{i<s<t}T(\widetilde{e_i},\widetilde{e_s},\widetilde{e_t})
c(\widetilde{e_i})c(\widetilde{e_s})c(\widetilde{e_t})
-\frac{1}{4}\sum^n_{i\neq s\neq t}T(\widetilde{e_i},\widetilde{e_s},\widetilde{e_t})
c(\widetilde{e_i})\widehat{c}(\widetilde{e_s})\widehat{c}(\widetilde{e_t})\Big) \nonumber\\
&+\Big(\frac{3}{2}\sum_{i<s<t}T(\widetilde{e_i},\widetilde{e_s},\widetilde{e_t})
c(\widetilde{e_i})c(\widetilde{e_s})c(\widetilde{e_t})
 -\frac{1}{4}\sum^n_{i\neq s\neq t}T(\widetilde{e_i},\widetilde{e_s},\widetilde{e_t})
c(\widetilde{e_i})\widehat{c}(\widetilde{e_s})\widehat{c}(\widetilde{e_t})\Big)c(\partial_{i})   \Big]\xi_{j}c(\xi).
\end{align}

  \begin{lem}
The following identities hold:
 \begin{align}
 &{\rm{Tr}}\Big(c(u)\widehat{c}(v)\widehat{c}(w)\sum^n_{i\neq s\neq t}T(\widetilde{e_i},\widetilde{e_s},\widetilde{e_t})
c(\widetilde{e_i})\widehat{c}(\widetilde{e_s})\widehat{c}(\widetilde{e_t})\Big)=-2T(u,v,w){\rm{Tr}}(\rm{Id});\\
&{\rm{Tr}}\Big(c(u)\widehat{c}(v)\widehat{c}(w)\sum_{i<s<t}T(\widetilde{e_i},\widetilde{e_s},\widetilde{e_t})
c(\widetilde{e_i})c(\widetilde{e_s})c(\widetilde{e_t})\Big)=0.
\end{align}
\end{lem}
\begin{proof}
By the relation of the Clifford action and $ {\rm{Tr}}(AB)= {\rm{Tr}}(BA) $, we have
 \begin{align}
 &{\rm{Tr}}\Big(c(u)\widehat{c}(v)\widehat{c}(w)\sum^n_{i\neq s\neq t}T(\widetilde{e_i},\widetilde{e_s},\widetilde{e_t})
c(\widetilde{e_i})\widehat{c}(\widetilde{e_s})\widehat{c}(\widetilde{e_t})\Big)\nonumber\\
= &\sum^n_{i\neq s\neq t}T(\widetilde{e_i},\widetilde{e_s},\widetilde{e_t}){\rm{Tr}}\Big(c(u)\widehat{c}(v)\widehat{c}(w)
c(\widetilde{e_i})\widehat{c}(\widetilde{e_s})\widehat{c}(\widetilde{e_t})\Big)\nonumber\\
= &-\sum^n_{i\neq s\neq t}T(\widetilde{e_i},\widetilde{e_s},\widetilde{e_t}){\rm{Tr}}\Big(c(u)\widehat{c}(v)\widehat{c}(w)
c(\widetilde{e_i})\widehat{c}(\widetilde{e_t})\widehat{c}(\widetilde{e_s})\Big)\nonumber\\
= &\sum^n_{i\neq s\neq t}T(\widetilde{e_i},\widetilde{e_s},\widetilde{e_t}){\rm{Tr}}\Big(c(u)\widehat{c}(v)\widehat{c}(w)
\widehat{c}(\widetilde{e_t})c(\widetilde{e_i})\widehat{c}(\widetilde{e_s})\Big)\nonumber\\
= &\sum^n_{i\neq s\neq t}T(\widetilde{e_i},\widetilde{e_s},\widetilde{e_t})
{\rm{Tr}}\Big(c(u)\widehat{c}(v)\big(-\widehat{c}(\widetilde{e_t})\widehat{c}(w)+2g(w,\widetilde{e_t})
\big)c(\widetilde{e_i})\widehat{c}(\widetilde{e_s})\Big)\nonumber\\
= &-\sum^n_{i\neq s\neq t}T(\widetilde{e_i},\widetilde{e_s},\widetilde{e_t})
{\rm{Tr}}\Big(c(u)\widehat{c}(v)\widehat{c}(\widetilde{e_t})\widehat{c}(w)
c(\widetilde{e_i})\widehat{c}(\widetilde{e_s})\Big)\nonumber\\
&+2g(w,\widetilde{e_t})\sum^n_{i\neq s\neq t}T(\widetilde{e_i},\widetilde{e_s},\widetilde{e_t})
{\rm{Tr}}\Big(c(u)\widehat{c}(v)c(\widetilde{e_i})\widehat{c}(\widetilde{e_s})\Big)\nonumber\\
=&\cdots\nonumber\\
= &-\sum^n_{i\neq s\neq t}T(\widetilde{e_i},\widetilde{e_s},\widetilde{e_t})
{\rm{Tr}}\Big(\widehat{c}(\widetilde{e_t})c(u)\widehat{c}(v)\widehat{c}(w)
c(\widetilde{e_i})\widehat{c}(\widetilde{e_s})\Big)\nonumber\\
&-2g(v,\widetilde{e_t})\sum^n_{i\neq s\neq t}T(\widetilde{e_i},\widetilde{e_s},\widetilde{e_t})
{\rm{Tr}}\Big(c(u)\widehat{c}(w)c(\widetilde{e_i})\widehat{c}(\widetilde{e_s})\Big)\nonumber\\
&+2g(w,\widetilde{e_t})\sum^n_{i\neq s\neq t}T(\widetilde{e_i},\widetilde{e_s},\widetilde{e_t})
{\rm{Tr}}\Big(c(u)\widehat{c}(v)c(\widetilde{e_i})\widehat{c}(\widetilde{e_s})\Big).
\end{align}
Then
 \begin{align}
 &{\rm{Tr}}\Big(c(u)\widehat{c}(v)\widehat{c}(w)\sum^n_{i\neq s\neq t}T(\widetilde{e_i},\widetilde{e_s},\widetilde{e_t})
c(\widetilde{e_i})\widehat{c}(\widetilde{e_s})\widehat{c}(\widetilde{e_t})\Big)\nonumber\\
=&-g(v,\widetilde{e_t})\sum^n_{i\neq s\neq t}T(\widetilde{e_i},\widetilde{e_s},\widetilde{e_t})
{\rm{Tr}}\Big(c(u)\widehat{c}(w)c(\widetilde{e_i})\widehat{c}(\widetilde{e_s})\Big)\nonumber\\
&+g(w,\widetilde{e_t})\sum^n_{i\neq s\neq t}T(\widetilde{e_i},\widetilde{e_s},\widetilde{e_t})
{\rm{Tr}}\Big(c(u)\widehat{c}(v)c(\widetilde{e_i})\widehat{c}(\widetilde{e_s})\Big)\nonumber\\
=&-2T(u,v,w){\rm{Tr}}(\rm{Id}).
\end{align}
On the other hand,
 \begin{align}
&{\rm{Tr}}\Big(c(u)\widehat{c}(v)\widehat{c}(w)\sum_{i<s<t}T(\widetilde{e_i},\widetilde{e_s},\widetilde{e_t})
c(\widetilde{e_i})c(\widetilde{e_s})c(\widetilde{e_t})\Big)\nonumber\\
=&\sum_{i<s<t}T(\widetilde{e_i},\widetilde{e_s},\widetilde{e_t}){\rm{Tr}}\Big(c(u)\widehat{c}(w)\widehat{c}(v)
c(\widetilde{e_i})c(\widetilde{e_s})c(\widetilde{e_t})\Big)\nonumber\\
=& \sum_{i<s<t}T(\widetilde{e_i},\widetilde{e_s},\widetilde{e_t}){\rm{Tr}}\Big(c(u)\big(-\widehat{c}(v)\widehat{c}(w)+2g(w,v)\big)
c(\widetilde{e_i})c(\widetilde{e_s})c(\widetilde{e_t})\Big)\nonumber\\
=& -\sum_{i<s<t}T(\widetilde{e_i},\widetilde{e_s},\widetilde{e_t}){\rm{Tr}}\Big(c(u)\widehat{c}(v)\widehat{c}(w)
c(\widetilde{e_i})c(\widetilde{e_s})c(\widetilde{e_t})\Big)\nonumber\\
&+ 2g(w,v)\sum_{i<s<t}T(\widetilde{e_i},\widetilde{e_s},\widetilde{e_t}){\rm{Tr}}\Big(c(u)
c(\widetilde{e_i})c(\widetilde{e_s})c(\widetilde{e_t})\Big).
\end{align}
By the generating relation, we see that the left hand side of (3.27) equals to
 \begin{align}
{\rm{Tr}}\Big(c(u)\widehat{c}(v)\widehat{c}(w)\sum_{i<s<t}T(\widetilde{e_i},\widetilde{e_s},\widetilde{e_t})
c(\widetilde{e_i})c(\widetilde{e_s})c(\widetilde{e_t})\Big)=0.
\end{align}
Then the proof of the lemma is complete.
\end{proof}
 \begin{lem}
The following identities hold:
 \begin{align}
 &\int_{|\xi|=1}{\rm{Tr}}\Big(\sum_{i,j}g^{ij}c(u)\widehat{c}(v)\widehat{c}(w)
 \sum_{l<s<t}T(\widetilde{e_l},\widetilde{e_s},\widetilde{e_t})c(\widetilde{e_l})c(\widetilde{e_s})c(\widetilde{e_t})
 c(\partial_{i})\xi_{j}c(\xi) \Big)\sigma(\xi)
 =0;\\
 &\int_{|\xi|=1}{\rm{Tr}}\Big(\sum_{i,j}g^{ij}c(u)\widehat{c}(v)\widehat{c}(w)c(\partial_{i})
 \sum_{l<s<t}T(\widetilde{e_l},\widetilde{e_s},\widetilde{e_t})c(\widetilde{e_l})c(\widetilde{e_s})c(\widetilde{e_t})
 \xi_{j}c(\xi) \Big)
 \sigma(\xi)
 =0.
\end{align}
\end{lem}
 \begin{proof}
By the relation of the Clifford action and $ {\rm{Tr}}(AB)= {\rm{Tr}}(BA) $, in terms of the condition (2.24), we get
\begin{align}
 &\int_{|\xi|=1}{\rm{Tr}}\Big(\sum_{i,j}g^{ij}c(u)\widehat{c}(v)\widehat{c}(w)
 \sum_{l<s<t}T(\widetilde{e_l},\widetilde{e_s},\widetilde{e_t})c(\widetilde{e_l})c(\widetilde{e_s})c(\widetilde{e_t})
 c(\partial_{i})\xi_{j}c(\xi) \Big)\sigma(\xi)\nonumber\\
=&\sum_{l<s<t}T(\widetilde{e_l},\widetilde{e_s},\widetilde{e_t})
\int_{|\xi|=1}{\rm{Tr}}\Big(\sum_{i,j}g^{ij}c(u)\widehat{c}(v)\widehat{c}(w)
c(\widetilde{e_l})c(\widetilde{e_s})c(\widetilde{e_t})
 c(\partial_{i})\xi_{j}c(\xi) \Big)\sigma(\xi)\nonumber\\
  =&\sum_{l<s<t}T(\widetilde{e_l},\widetilde{e_s},\widetilde{e_t})\sum_{i,h}\delta _{i}^{h}\frac{1}{n}{\rm{vol}}(S^{n-1})
  {\rm{Tr}}\Big(c(u)\widehat{c}(v)\widehat{c}(w)
c(\widetilde{e_l})c(\widetilde{e_s})c(\widetilde{e_t})
 c(\partial_{i}) c(\widetilde{e_h}) \Big)\nonumber\\
   =&-\sum_{l<s<t}T(\widetilde{e_l},\widetilde{e_s},\widetilde{e_t})\sum_{i}\frac{1}{n}{\rm{vol}}(S^{n-1})
  {\rm{Tr}}\Big(c(u)\widehat{c}(v)\widehat{c}(w)
c(\widetilde{e_l})c(\widetilde{e_s})c(\widetilde{e_t})
  \Big)\nonumber\\
   =&-\sum_{l<s<t}T(\widetilde{e_l},\widetilde{e_s},\widetilde{e_t})\sum_{i}\frac{1}{n}{\rm{vol}}(S^{n-1})
  {\rm{Tr}}\Big(c(u)\big(-\widehat{c}(w)\widehat{c}(v)+2g(v,w)\big)
c(\widetilde{e_l})c(\widetilde{e_s})c(\widetilde{e_t})
  \Big).
\end{align}
Then
  \begin{align}
 &\int_{|\xi|=1}{\rm{Tr}}\Big(\sum_{i,j}g^{ij}c(u)\widehat{c}(v)\widehat{c}(w)
 \sum_{l<s<t}T(\widetilde{e_l},\widetilde{e_s},\widetilde{e_t})c(\widetilde{e_l})c(\widetilde{e_s})c(\widetilde{e_t})
 c(\partial_{i})\xi_{j}c(\xi) \Big)\sigma(\xi)\nonumber\\
=&-g(v,w)\sum_{l<s<t}T(\widetilde{e_l},\widetilde{e_s},\widetilde{e_t})\sum_{i}\frac{1}{n}{\rm{vol}}(S^{n-1})
  {\rm{Tr}}\Big(c(u)c(\widetilde{e_l})c(\widetilde{e_s})c(\widetilde{e_t})
  \Big)\nonumber\\
  =&0,
\end{align}
and
  \begin{align}
 &\int_{|\xi|=1}{\rm{Tr}}\Big(\sum_{i,j}g^{ij}c(u)\widehat{c}(v)\widehat{c}(w)
 c(\partial_{i}) \sum_{l<s<t}T(\widetilde{e_l},\widetilde{e_s},\widetilde{e_t})c(\widetilde{e_l})c(\widetilde{e_s})c(\widetilde{e_t})
\xi_{j}c(\xi) \Big)\sigma(\xi)\nonumber\\
  =&\sum_{l<s<t}T(\widetilde{e_l},\widetilde{e_s},\widetilde{e_t})\sum_{i,h}\delta _{i}^{h}\frac{1}{n}{\rm{vol}}(S^{n-1})
  {\rm{Tr}}\Big(c(u)\widehat{c}(v)\widehat{c}(w) c(\partial_{i})
c(\widetilde{e_l})c(\widetilde{e_s})c(\widetilde{e_t})
c(\widetilde{e_h}) \Big)\nonumber\\
   = &\sum_{l<s<t}T(\widetilde{e_l},\widetilde{e_s},\widetilde{e_t})\sum_{i } \frac{1}{n}{\rm{vol}}(S^{n-1})
  {\rm{Tr}}\Big(\big(-c(\widetilde{e_i})c(u)-2g(u,\widetilde{e_i})\big)\widehat{c}(v)\widehat{c}(w)
c(\widetilde{e_l})c(\widetilde{e_s})c(\widetilde{e_t})
c(\widetilde{e_i}) \Big)\nonumber\\
  = &\sum_{l<s<t}T(\widetilde{e_l},\widetilde{e_s},\widetilde{e_t})\sum_{i } \frac{1}{n}{\rm{vol}}(S^{n-1})
  {\rm{Tr}}\Big( c(u) \widehat{c}(v)\widehat{c}(w)
c(\widetilde{e_l})c(\widetilde{e_s})c(\widetilde{e_t})
  \Big)\nonumber\\
 &-2g(u,\widetilde{e_i})\sum_{l<s<t}T(\widetilde{e_l},\widetilde{e_s},\widetilde{e_t})\sum_{i } \frac{1}{n}{\rm{vol}}(S^{n-1})
  {\rm{Tr}}\Big(\widehat{c}(v)\widehat{c}(w)c(\widetilde{e_l})c(\widetilde{e_s})c(\widetilde{e_t})
c(\widetilde{e_i}) \Big)\nonumber\\
= &0.
\end{align}
 \end{proof}

 \begin{lem}
The following identities hold:
 \begin{align}
 &\int_{|\xi|=1}{\rm{Tr}}\Big(\sum_{i,j}g^{ij}c(u)\widehat{c}(v)\widehat{c}(w)c(\partial_{i})\sum_{l\neq s\neq t}
 T(\widetilde{e_l},\widetilde{e_s},\widetilde{e_t})
 c(\widetilde{e_l})\widehat{c}(\widetilde{e_s})\widehat{c}(\widetilde{e_t})\xi_{j}
 c(\xi) \Big)\sigma(\xi)\nonumber\\
 =&2T(u,v,w){\rm{Tr}}(\rm{Id}){\rm{vol}}(S^{n-1});\\
 &\int_{|\xi|=1}{\rm{Tr}}\Big(\sum_{i,j}g^{ij}c(u)\widehat{c}(v)\widehat{c}(w)\sum_{l\neq s\neq t}T(\widetilde{e_l},\widetilde{e_s},\widetilde{e_t})
 c(\widetilde{e_l})\widehat{c}(\widetilde{e_s})\widehat{c}(\widetilde{e_t})c(\partial_{i})\xi_{j}c(\xi) \Big)\sigma(\xi)
\nonumber\\
 =&2T(u,v,w){\rm{Tr}}(\rm{Id}){\rm{vol}}(S^{n-1}).
\end{align}
\end{lem}
\begin{proof}
By the relation of the Clifford action and $ {\rm{Tr}}(AB)= {\rm{Tr}}(BA) $, in terms of the condition (2.24), we get
 \begin{align}
 &\int_{|\xi|=1}{\rm{Tr}}\Big(\sum_{i,j}g^{ij}c(u)\widehat{c}(v)\widehat{c}(w)c(\partial_{i})\sum_{l\neq s\neq t}
 T(\widetilde{e_l},\widetilde{e_s},\widetilde{e_t})
 c(\widetilde{e_l})\widehat{c}(\widetilde{e_s})\widehat{c}(\widetilde{e_t})\xi_{j}
 c(\xi) \Big)\sigma(\xi)\nonumber\\
 =&\int_{|\xi|=1}\sum_{l\neq s\neq t}
 T(\widetilde{e_l},\widetilde{e_s},\widetilde{e_t}){\rm{Tr}}\Big(\sum_{i } c(u)\widehat{c}(v)\widehat{c}(w)c(\partial_{i})
 c(\widetilde{e_l})\widehat{c}(\widetilde{e_s})\widehat{c}(\widetilde{e_t})\xi_{i}
 c(\xi) \Big)\sigma(\xi)\nonumber\\
  =&\sum_{i,h}\delta _{i}^{h}\frac{1}{n}{\rm{vol}}(S^{n-1})\sum_{l\neq s\neq t}
 T(\widetilde{e_l},\widetilde{e_s},\widetilde{e_t}){\rm{Tr}}\Big(\sum_{i } c(u)\widehat{c}(v)\widehat{c}(w)c(\partial_{i})
 c(\widetilde{e_l})\widehat{c}(\widetilde{e_s})\widehat{c}(\widetilde{e_t})
 c(\widetilde{e_h}) \Big) \nonumber\\
  =&\sum_{i } \frac{1}{n}{\rm{vol}}(S^{n-1})\sum_{l\neq s\neq t}
 T(\widetilde{e_l},\widetilde{e_s},\widetilde{e_t}){\rm{Tr}}\Big(\sum_{i } c(u)\widehat{c}(v)\widehat{c}(w)c(\partial_{i})
 c(\widetilde{e_l})\widehat{c}(\widetilde{e_s})\widehat{c}(\widetilde{e_t})
 c(\widetilde{e_i}) \Big) \nonumber\\
   =&\sum_{i } \frac{1}{n}{\rm{vol}}(S^{n-1})\sum_{l\neq s\neq t}
 T(\widetilde{e_l},\widetilde{e_s},\widetilde{e_t}){\rm{Tr}}\Big(\sum_{i } c(u)\widehat{c}(v)\widehat{c}(w)c(\partial_{i})
 \big(-c(\widetilde{e_i})c(\widetilde{e_l})-2\delta_{l}^{i}\big)\widehat{c}(\widetilde{e_s})\widehat{c}(\widetilde{e_t})
\Big) \nonumber\\
=&\sum_{i } \frac{1}{n}{\rm{vol}}(S^{n-1})\sum_{l\neq s\neq t}
 T(\widetilde{e_l},\widetilde{e_s},\widetilde{e_t}){\rm{Tr}}\Big(\sum_{i } c(u)\widehat{c}(v)\widehat{c}(w)
 c(\widetilde{e_l}) \widehat{c}(\widetilde{e_s})\widehat{c}(\widetilde{e_t})
\Big) \nonumber\\
&-2\sum_{i }\delta_{l}^{i} \frac{1}{n}{\rm{vol}}(S^{n-1})\sum_{l\neq s\neq t}
 T(\widetilde{e_l},\widetilde{e_s},\widetilde{e_t}){\rm{Tr}}\Big(\sum_{i } c(u)\widehat{c}(v)\widehat{c}(w)c(\partial_{i})
 \widehat{c}(\widetilde{e_s})\widehat{c}(\widetilde{e_t})
\Big) \nonumber\\
=&-\sum_{i } \frac{1}{n}{\rm{vol}}(S^{n-1})\sum_{l\neq s\neq t}
 T(\widetilde{e_l},\widetilde{e_s},\widetilde{e_t}){\rm{Tr}}\Big(\sum_{i } c(u)\widehat{c}(v)\widehat{c}(w)
 c(\widetilde{e_l}) \widehat{c}(\widetilde{e_s})\widehat{c}(\widetilde{e_t})
\Big) \nonumber\\
=&2T(u,v,w){\rm{Tr}}(\rm{Id}){\rm{vol}}(S^{n-1}).
 \end{align}
In the same way we obtain
 \begin{align}
  &\int_{|\xi|=1}{\rm{Tr}}\Big(\sum_{i,j}g^{ij}c(u)\widehat{c}(v)\widehat{c}(w)\sum_{l\neq s\neq t}T(\widetilde{e_l},\widetilde{e_s},\widetilde{e_t})
 c(\widetilde{e_l})\widehat{c}(\widetilde{e_s})\widehat{c}(\widetilde{e_t})c(\partial_{i})\xi_{j}c(\xi) \Big)\sigma(\xi)
\nonumber\\
=&\sum_{l\neq s\neq t}T(\widetilde{e_l},\widetilde{e_s},\widetilde{e_t})\int_{|\xi|=1}{\rm{Tr}}\Big(\sum_{i } c(u)\widehat{c}(v)\widehat{c}(w)
 c(\widetilde{e_l})\widehat{c}(\widetilde{e_s})\widehat{c}(\widetilde{e_t})c(\partial_{i})\xi_{i}c(\xi) \Big)\sigma(\xi)
\nonumber\\
=&\sum_{i,h}\delta _{i}^{h}\frac{1}{n}{\rm{vol}}(S^{n-1})\sum_{l\neq s\neq t}T(\widetilde{e_l},\widetilde{e_s},\widetilde{e_t})
  {\rm{Tr}}\Big( c(u)\widehat{c}(v)\widehat{c}(w)
 c(\widetilde{e_l})\widehat{c}(\widetilde{e_s})\widehat{c}(\widetilde{e_t})c(\partial_{i})c(\widetilde{e_h}) \Big)\nonumber\\
 =&  \sum_{i}\frac{1}{n}{\rm{vol}}(S^{n-1})\sum_{l\neq s\neq t}T(\widetilde{e_l},\widetilde{e_s},\widetilde{e_t})
  {\rm{Tr}}\Big( c(u)\widehat{c}(v)\widehat{c}(w)
 c(\widetilde{e_l})\widehat{c}(\widetilde{e_s})\widehat{c}(\widetilde{e_t}) \Big)\nonumber\\
 =&2T(u,v,w){\rm{Tr}}(\rm{Id}){\rm{vol}}(S^{n-1}).
\end{align}
It follows from Lemma 3.7, Lemma 3.8 and Lemma 3.9, we obtain
\begin{align}
&\mathscr{T}_{3}(c(u),\widehat{c}(v),\widehat{c}(w))\nonumber\\
&\int_{|\xi|=1}{\rm Tr}_{\wedge^{*}(T^{*}M)}
[\sigma_{-2m}\big(c(u)\widehat{c}(v)\widehat{c}(w)(d+\delta)_{T }^{-2m+1}\big)](x_{0})|_{|\xi|=1}\sigma(\xi) \nonumber\\
=&\int_{|\xi|=1}{\rm Tr}_{\wedge^{*}(T^{*}M)}\Big[
c(u)\widehat{c}(v)\widehat{c}(w)\Big(-\frac{3}{2}\sum_{i<s<t}T(\widetilde{e_i},\widetilde{e_s},\widetilde{e_t})
c(\widetilde{e_i})c(\widetilde{e_s})c(\widetilde{e_t})
 +\frac{1}{4}\sum^n_{i\neq s\neq t}T(\widetilde{e_i},\widetilde{e_s},\widetilde{e_t})
c(\widetilde{e_i})\widehat{c}(\widetilde{e_s})\widehat{c}(\widetilde{e_t})\Big)\nonumber\\
&+m\sum_{i,j}g^{ij}c(u)\widehat{c}(v)\widehat{c}(w)\Big[c(\partial_{i})\Big(-\frac{3}{2}\sum_{i<s<t}T(\widetilde{e_i},\widetilde{e_s},\widetilde{e_t})
c(\widetilde{e_i})c(\widetilde{e_s})c(\widetilde{e_t})
+\frac{1}{4}\sum^n_{i\neq s\neq t}T(\widetilde{e_i},\widetilde{e_s},\widetilde{e_t})
c(\widetilde{e_i})\widehat{c}(\widetilde{e_s})\widehat{c}(\widetilde{e_t})\Big) \nonumber\\
&+\Big(-\frac{3}{2}\sum_{i<s<t}T(\widetilde{e_i},\widetilde{e_s},\widetilde{e_t})
c(\widetilde{e_i})c(\widetilde{e_s})c(\widetilde{e_t})
+\frac{1}{4}\sum^n_{i\neq s\neq t}T(\widetilde{e_i},\widetilde{e_s},\widetilde{e_t})
c(\widetilde{e_i})\widehat{c}(\widetilde{e_s})\widehat{c}(\widetilde{e_t})\Big)c(\partial_{i})   \Big]\xi_{j}c(\xi)\Big](x_{0})|_{|\xi|=1}\sigma(\xi)
 \nonumber\\
 =&2^{2m-1}(1-2m)T(u,v,w)
{\rm{vol}}(S^{2m-1}).
\end{align}
\end{proof}
Substituting (3.38) into (3.7) leads to the desired equality, and the proof of Theorem 1.3 is complete.
\begin{defn}
Let $Q$ be a first order differential operator acting on $\mathcal{H}$, and $Q^{*}$ be the self adjoint operator of $Q$.
Let $D=Q+Q^{*}$ and $\widetilde{D}=\sqrt{-1}(Q-Q^{*})$ are the elliptic operators and $a \in \mathcal{A}$, where $\mathcal{A}$
is a $C^{*}$ alegbra,
and  $[D,a] $, $[\widetilde{D},a]$ are bounded operators.
Then we obtain the de-Rham Hodge  spectral triple $(\mathcal{A},\mathcal{H},\widetilde{D})$,
  the spectral torsion can be defined by $Wres[[D,a][\widetilde{D},a^{2}][\widetilde{D},a^{3}]D^{-2m+1}]$.
\end{defn}
\begin{exam} Let $L=\frac{3}{2}\sum_{i<s<t}T(\widetilde{e_i},\widetilde{e_s},\widetilde{e_t})
c(\widetilde{e_i})c(\widetilde{e_s})c(\widetilde{e_t})
 -\frac{1}{4}\sum^n_{i\neq s\neq t}T(\widetilde{e_i},\widetilde{e_s},\widetilde{e_t})
c(\widetilde{e_i})\widehat{c}(\widetilde{e_s})\widehat{c}(\widetilde{e_t})$,
then $L^{*}=L$. Let $Q=d+\frac{L}{2} $, $Q^{*}=\delta+\frac{L}{2} $,  $D=Q+Q^{*}=(d+\delta)+L$ and
$\widetilde{D}=\sqrt{-1}(Q-Q^{*})=\sqrt{-1}(d-\delta)$,
then we get  $[D,f]=c(df)$ and  $[\widetilde{D},f]=\sqrt{-1}\widehat{c}(df)$.
\end{exam}

\section{The four-linear functional  for the  de-Rham Hodge type operator}

In this section we want to consider the four-linear functional  for the  de-Rham Hodge type operator with torsion.
To avoid technique terminology we only state our results for
compact oriented Riemannian manifold of even dimension $n=2m$ by using the trace of  de-Rham Hodge type operator and the noncommutative residue density.
Using an explicit formula for $d+\delta$, we can reformulate $d+\delta +\sqrt{-1}T_{3}$  and $d+\delta +\sqrt{-1}T_{4}$  as follows.
\begin{defn}
The de-Rham Hodge type operator  with 4-form pertubation $d+\delta +\sqrt{-1}T_{3}$ is the first order differential operator
on $ \Gamma(M,\wedge^{*}(T^{*}M))$  given by the formula
 \begin{align}
d+\delta +\sqrt{-1}T_{3}=&d+\delta
+\sum_{k<l <\alpha<\beta}\sqrt{-1}T(\widetilde{e_k},\widetilde{e_l},\widetilde{e_\alpha},\widetilde{e_\beta})
c(\widetilde{e_k})c(\widetilde{e_l})\widehat{c}(\widetilde{e_\alpha})\widehat{c}(\widetilde{e_\beta}).
\end{align}
\end{defn}
\begin{defn}
The de-Rham Hodge type operator  with 4-form pertubation $d+\delta +T_{4}$ is the first order differential operator
on $ \Gamma(M,\wedge^{*}(T^{*}M))$  given by the formula
 \begin{align}
d+\delta +\sqrt{-1}T_{4}=&d+\delta
+\sum_{k<l <\alpha<\beta}\sqrt{-1}T(\widetilde{e_k},\widetilde{e_l},\widetilde{e_\alpha},\widetilde{e_\beta})
\widehat{c}(\widetilde{e_k})\widehat{c}(\widetilde{e_l})\widehat{c}(\widetilde{e_\alpha})\widehat{c}(\widetilde{e_\beta}).
\end{align}
\end{defn}

\begin{defn}\cite{DSZ}
For $d+\delta +\sqrt{-1}T_{3}$, the four-linear Clifford multiplication by functional of differential one-forms
$c(u),c(v),\widehat{c}(w),\widehat{c}(z)$
\begin{align}
\mathscr{T}_{4}(c(u),c(v),\widehat{c}(w),\widehat{c}(z))={\rm Wres}\big(c(u)c(v)\widehat{c}(w)\widehat{c}(z)(d+\delta+\sqrt{-1}T_{3})^{-n+1}\big)
\end{align}
is called the spectral four form $\mathscr{T}_{4}$ associated with the de-Rham Hodge type operator.
\end{defn}

\begin{defn}\cite{DSZ}
For $d+\delta +\sqrt{-1}T_{4}$, the four-linear Clifford multiplication by functional of differential one-forms
$\widehat{c}(u),\widehat{c}(v),\widehat{c}(w),\widehat{c}(z)$
\begin{align}
\mathscr{T}_{5}(\widehat{c}(u),\widehat{c}(v),\widehat{c}(w),\widehat{c}(z))
={\rm Wres}\big(\widehat{c}(u)\widehat{c}(v)\widehat{c}(w)\widehat{c}(z)(d+\delta+\sqrt{-1}T_{4})^{-n+1}\big)
\end{align}
is called the spectral four form  $\mathscr{T}_{5}$ associated with the de-Rham Hodge type operator.
\end{defn}

\subsection{The proof of Theorem 1.4}
For any fixed point $x_0\in M$, we can choose the normal coordinates
$U$ of $x_0$ in $M$ and compute $\sigma_{-n}\big(c(u)c(v)c(w)(d+\delta)_{T_{3}}^{-n+1}\big)$.
Let $c(\widetilde{e_j})=\epsilon (\widetilde{e_j*})-\iota (\widetilde{e_j*});\
\widehat{c}(\widetilde{e_j})=\epsilon (\widetilde{e_j*})+\iota (\widetilde{e_j*})$ act on $\Gamma(M,\wedge^{*}(T^{*}M))$.
Then we have $ \partial_{x_j}[c(\widetilde{e_{i}})]=\partial_{x_j}[\widehat{c}(\widetilde{e_{i}})]=0$ in the above frame
$\{\widetilde{e_{1}},\widetilde{e_{2}},\cdots,\widetilde{e_{n}}\}$. In terms of normal coordinates about $x_{0}$ one has:
$\sigma^{j}_{S(TM)}(x_{0})=0$ $e_{j}\big(c(e_{i})\big)(x_{0})=0$, $\Gamma^{k}(x_{0})=0$.
Substituting above results into the formula (2.18), we obtain
\begin{align}
 &\sigma_{-2m}(c(u)c(v)\widehat{c}(w)\widehat{c}(z)(d+\delta)_{T_{3}}^{1-2m})(x_{0})|_{|\xi|=1} \nonumber\\
 =&c(u)c(v)\widehat{c}(w)\widehat{c}(z)\Big\{|\xi|^{-2m}\sigma_0((d+\delta)_{T_{3}})+\sum_{j=1}^{2m+2}\partial_{\xi_j}(|\xi|^{-2m})\partial_{x_j}c(\xi)+
\Big[m\sigma_2((d+\delta)_{T_{3}})^{-m+1}\sigma_{-3}((d+\delta)_{T_{3}}^{-2})  \nonumber\\
&-\sqrt{-1} \sum_{k=0}^{m-2}\sum_{\mu=1}^{2m+2}
\partial_{\xi_{\mu}}\sigma_{2}^{-m+k+1}((d+\delta)_{T_{3}}^2)
\partial_{x_{\mu}}\sigma_{2}^{-1}((d+\delta)_{T_{3}}^2)(\sigma_2((d+\delta)_{T_{3}}^2))^{-k}\Big ]\sqrt{-1}c(\xi)\Big\}(x_{0})|_{|\xi|=1} \nonumber\\
=&c(u)c(v)\widehat{c}(w)\widehat{c}(z)\sqrt{-1}\sum_{k<l <\alpha<\beta}T(\widetilde{e_k},\widetilde{e_l},\widetilde{e_\alpha},\widetilde{e_\beta})
c(\widetilde{e_k})c(\widetilde{e_l})\widehat{c}(\widetilde{e_\alpha})\widehat{c}(\widetilde{e_\beta}) \nonumber\\
&+m\sum_{i,j}g^{ij}c(u)c(v)\widehat{c}(w)\widehat{c}(z)\sqrt{-1}\Big[c(\partial_{i}) \sum_{k<l <\alpha<\beta}T(\widetilde{e_k},\widetilde{e_l},\widetilde{e_\alpha},\widetilde{e_\beta})
c(\widetilde{e_k})c(\widetilde{e_l})\widehat{c}(\widetilde{e_\alpha})\widehat{c}(\widetilde{e_\beta})
 \nonumber\\
&+ \sum_{k<l <\alpha<\beta}T(\widetilde{e_k},\widetilde{e_l},\widetilde{e_\alpha},\widetilde{e_\beta})
c(\widetilde{e_k})c(\widetilde{e_l})\widehat{c}(\widetilde{e_\alpha})\widehat{c}(\widetilde{e_\beta}) c(\partial_{i})   \Big]\xi_{j}c(\xi).
\end{align}

Let $u=\sum_{i=1}^{n}u_{i}\widetilde{e_{i}}$, $v=\sum_{j=1}^{n}v_{j}\widetilde{e_{j}} $, $w=\sum_{k=1}^{n}w_{k}\widetilde{e_{k}} $,
$z=\sum_{l=1}^{n}z_{l}\widetilde{e_{l}} $,
where $\{\widetilde{e_{1}},\widetilde{e_{2}},\cdots,\widetilde{e_{n}}\}$ is the orthogonal basis about $g^{TM}$,
then  $c(u)= \sum_{i=1}^{n}u_{i}c(\widetilde{e_{i}})$,$c(v)= \sum_{j=1}^{n}v_{j}c(\widetilde{e_{j}})$,
$c(w)= \sum_{k=1}^{n}w_{k}c(\widetilde{e_{k}})$, $c(z)= \sum_{l=1}^{n}z_{l}c(\widetilde{e_{l}})$.
The following Lemmas of traces in terms of the Clifford action
is very efficient for solving the spectral torsion  associated with the de-Rham Hodge type operator.
 \begin{lem}
The following identities hold:
 \begin{align}
 &{\rm{Tr}}\Big(\sum_{k<l <\alpha<\beta}T(\widetilde{e_k},\widetilde{e_l},\widetilde{e_\alpha},\widetilde{e_\beta})c(u)c(v)\widehat{c}(w)\widehat{c}(z)
c(\widetilde{e_k})c(\widetilde{e_l})\widehat{c}(\widetilde{e_\alpha})\widehat{c}(\widetilde{e_\beta})\Big)(x_{0})  \nonumber\\
=&-\frac{1}{6}T(u,v,w,z){\rm{Tr}}(\rm{Id}).
\end{align}
\end{lem}
\begin{proof}
By the relation of the Clifford action and $ {\rm{Tr}}(AB)= {\rm{Tr}}(BA) $, we get
 \begin{align}
 &{\rm{Tr}}\Big(\sum_{k<l <\alpha<\beta}T(\widetilde{e_k},\widetilde{e_l},\widetilde{e_\alpha},\widetilde{e_\beta})c(u)c(v)\widehat{c}(w)\widehat{c}(z)
c(\widetilde{e_k})c(\widetilde{e_l})\widehat{c}(\widetilde{e_\alpha})\widehat{c}(\widetilde{e_\beta})\Big)(x_{0}) \nonumber\\
= &\sum_{k<l <\alpha<\beta}T(\widetilde{e_k},\widetilde{e_l},\widetilde{e_\alpha},\widetilde{e_\beta})\sum_{a,b,c,d }u_{a}v_{b}w_{c}z_{d}
{\rm{Tr}}\Big(c(\widetilde{e_a})c(\widetilde{e_b})\widehat{c}(\widetilde{e_c})\widehat{c}(\widetilde{e_d})
c(\widetilde{e_k})c(\widetilde{e_l})\widehat{c}(\widetilde{e_\alpha})\widehat{c}(\widetilde{e_\beta})\Big)(x_{0}) \nonumber\\
= &\sum_{k<l <\alpha<\beta}T(\widetilde{e_k},\widetilde{e_l},\widetilde{e_\alpha},\widetilde{e_\beta})\sum_{a,b,c,d }u_{a}v_{b}w_{c}z_{d}
{\rm{Tr}}\Big[c(\widetilde{e_a})\Big(-c(\widetilde{e_k})c(\widetilde{e_b})-2\delta_{k}^{b}   \Big)
\widehat{c}(\widetilde{e_c})\widehat{c}(\widetilde{e_d})
c(\widetilde{e_l})\widehat{c}(\widetilde{e_\alpha})\widehat{c}(\widetilde{e_\beta})\Big](x_{0}) \nonumber\\
= &-\sum_{k<l <\alpha<\beta}T(\widetilde{e_k},\widetilde{e_l},\widetilde{e_\alpha},\widetilde{e_\beta})\sum_{a,b,c,d }u_{a}v_{b}w_{c}z_{d}
{\rm{Tr}}\Big[c(\widetilde{e_a}) c(\widetilde{e_k})c(\widetilde{e_b})
\widehat{c}(\widetilde{e_c})\widehat{c}(\widetilde{e_d})
c(\widetilde{e_l})\widehat{c}(\widetilde{e_\alpha})\widehat{c}(\widetilde{e_\beta})\Big](x_{0}) \nonumber\\
&-\sum_{k<l <\alpha<\beta}T(\widetilde{e_k},\widetilde{e_l},\widetilde{e_\alpha},\widetilde{e_\beta})\sum_{a,b,c,d }u_{a}v_{b}w_{c}z_{d}
2\delta_{k}^{b}{\rm{Tr}}\Big[c(\widetilde{e_a})
\widehat{c}(\widetilde{e_c})\widehat{c}(\widetilde{e_d})
c(\widetilde{e_l})\widehat{c}(\widetilde{e_\alpha})\widehat{c}(\widetilde{e_\beta})\Big](x_{0}) \nonumber\\
= &-\sum_{k<l <\alpha<\beta}T(\widetilde{e_k},\widetilde{e_l},\widetilde{e_\alpha},\widetilde{e_\beta})\sum_{a,b,c,d }u_{a}v_{b}w_{c}z_{d}
{\rm{Tr}}\Big[\Big(-c(\widetilde{e_k})c(\widetilde{e_a})-2\delta_{k}^{a}\Big)c(\widetilde{e_b})
\widehat{c}(\widetilde{e_c})\widehat{c}(\widetilde{e_d})
c(\widetilde{e_l})\widehat{c}(\widetilde{e_\alpha})\widehat{c}(\widetilde{e_\beta})\Big](x_{0}) \nonumber\\
&-\sum_{k<l <\alpha<\beta}T(\widetilde{e_k},\widetilde{e_l},\widetilde{e_\alpha},\widetilde{e_\beta})\sum_{a,b,c,d }u_{a}v_{b}w_{c}z_{d}
2\delta_{k}^{b}{\rm{Tr}}\Big[c(\widetilde{e_a})
\widehat{c}(\widetilde{e_c})\widehat{c}(\widetilde{e_d})
c(\widetilde{e_l})\widehat{c}(\widetilde{e_\alpha})\widehat{c}(\widetilde{e_\beta})\Big](x_{0}) \nonumber\\
= &\sum_{k<l <\alpha<\beta}T(\widetilde{e_k},\widetilde{e_l},\widetilde{e_\alpha},\widetilde{e_\beta})\sum_{a,b,c,d }u_{a}v_{b}w_{c}z_{d}
{\rm{Tr}}\Big[c(\widetilde{e_k})c(\widetilde{e_a})c(\widetilde{e_b})
\widehat{c}(\widetilde{e_c})\widehat{c}(\widetilde{e_d})
c(\widetilde{e_l})\widehat{c}(\widetilde{e_\alpha})\widehat{c}(\widetilde{e_\beta})\Big](x_{0}) \nonumber\\
&+\sum_{k<l <\alpha<\beta}T(\widetilde{e_k},\widetilde{e_l},\widetilde{e_\alpha},\widetilde{e_\beta})\sum_{a,b,c,d }u_{a}v_{b}w_{c}z_{d}
2\delta_{k}^{a}{\rm{Tr}}\Big[c(\widetilde{e_b})
\widehat{c}(\widetilde{e_c})\widehat{c}(\widetilde{e_d})
c(\widetilde{e_l})\widehat{c}(\widetilde{e_\alpha})\widehat{c}(\widetilde{e_\beta})\Big](x_{0}) \nonumber\\
&-\sum_{k<l <\alpha<\beta}T(\widetilde{e_k},\widetilde{e_l},\widetilde{e_\alpha},\widetilde{e_\beta})\sum_{a,b,c,d }u_{a}v_{b}w_{c}z_{d}
2\delta_{k}^{b}{\rm{Tr}}\Big[c(\widetilde{e_a})
\widehat{c}(\widetilde{e_c})\widehat{c}(\widetilde{e_d})
c(\widetilde{e_l})\widehat{c}(\widetilde{e_\alpha})\widehat{c}(\widetilde{e_\beta})\Big](x_{0}) \nonumber\\
= &\sum_{k<l <\alpha<\beta}T(\widetilde{e_k},\widetilde{e_l},\widetilde{e_\alpha},\widetilde{e_\beta})\sum_{a,b,c,d }u_{a}v_{b}w_{c}z_{d}
{\rm{Tr}}\Big[c(\widetilde{e_a})c(\widetilde{e_b})
\widehat{c}(\widetilde{e_c})\widehat{c}(\widetilde{e_d})
\Big(-c(\widetilde{e_k})c(\widetilde{e_l}) -2\delta_{k}^{l}\Big)\widehat{c}(\widetilde{e_\alpha})\widehat{c}(\widetilde{e_\beta})\Big](x_{0}) \nonumber\\
&+\sum_{k<l <\alpha<\beta}T(\widetilde{e_k},\widetilde{e_l},\widetilde{e_\alpha},\widetilde{e_\beta})\sum_{a,b,c,d }u_{a}v_{b}w_{c}z_{d}
2\delta_{k}^{a}{\rm{Tr}}\Big[c(\widetilde{e_b})
\widehat{c}(\widetilde{e_c})\widehat{c}(\widetilde{e_d})
c(\widetilde{e_l})\widehat{c}(\widetilde{e_\alpha})\widehat{c}(\widetilde{e_\beta})\Big](x_{0}) \nonumber\\
&-\sum_{k<l <\alpha<\beta}T(\widetilde{e_k},\widetilde{e_l},\widetilde{e_\alpha},\widetilde{e_\beta})\sum_{a,b,c,d }u_{a}v_{b}w_{c}z_{d}
2\delta_{k}^{b}{\rm{Tr}}\Big[c(\widetilde{e_a})
\widehat{c}(\widetilde{e_c})\widehat{c}(\widetilde{e_d})
c(\widetilde{e_l})\widehat{c}(\widetilde{e_\alpha})\widehat{c}(\widetilde{e_\beta})\Big](x_{0}) \nonumber\\
= &-\sum_{k<l <\alpha<\beta}T(\widetilde{e_k},\widetilde{e_l},\widetilde{e_\alpha},\widetilde{e_\beta})\sum_{a,b,c,d }u_{a}v_{b}w_{c}z_{d}
{\rm{Tr}}\Big[c(\widetilde{e_a})c(\widetilde{e_b})
\widehat{c}(\widetilde{e_c})\widehat{c}(\widetilde{e_d})
 c(\widetilde{e_k})c(\widetilde{e_l})  \widehat{c}(\widetilde{e_\alpha})\widehat{c}(\widetilde{e_\beta})\Big](x_{0}) \nonumber\\
&-\sum_{k<l <\alpha<\beta}T(\widetilde{e_k},\widetilde{e_l},\widetilde{e_\alpha},\widetilde{e_\beta})\sum_{a,b,c,d }u_{a}v_{b}w_{c}z_{d}
2\delta_{k}^{l}{\rm{Tr}}\Big[c(\widetilde{e_a})c(\widetilde{e_b})
\widehat{c}(\widetilde{e_c})\widehat{c}(\widetilde{e_d})
 \widehat{c}(\widetilde{e_\alpha})\widehat{c}(\widetilde{e_\beta})\Big](x_{0}) \nonumber\\
&+\sum_{k<l <\alpha<\beta}T(\widetilde{e_k},\widetilde{e_l},\widetilde{e_\alpha},\widetilde{e_\beta})\sum_{a,b,c,d }u_{a}v_{b}w_{c}z_{d}
2\delta_{k}^{a}{\rm{Tr}}\Big[c(\widetilde{e_b})
\widehat{c}(\widetilde{e_c})\widehat{c}(\widetilde{e_d})
c(\widetilde{e_l})\widehat{c}(\widetilde{e_\alpha})\widehat{c}(\widetilde{e_\beta})\Big](x_{0}) \nonumber\\
&-\sum_{k<l <\alpha<\beta}T(\widetilde{e_k},\widetilde{e_l},\widetilde{e_\alpha},\widetilde{e_\beta})\sum_{a,b,c,d }u_{a}v_{b}w_{c}z_{d}
2\delta_{k}^{b}{\rm{Tr}}\Big[c(\widetilde{e_a})
\widehat{c}(\widetilde{e_c})\widehat{c}(\widetilde{e_d})
c(\widetilde{e_l})\widehat{c}(\widetilde{e_\alpha})\widehat{c}(\widetilde{e_\beta})\Big](x_{0}).
\end{align}
By the generating relation, we see that the left hand side of (4.7) equals
 \begin{align}
 &{\rm{Tr}}\Big(\sum_{k<l <\alpha<\beta}T(\widetilde{e_k},\widetilde{e_l},\widetilde{e_\alpha},\widetilde{e_\beta})c(u)c(v)\widehat{c}(w)\widehat{c}(z)
c(\widetilde{e_k})c(\widetilde{e_l})\widehat{c}(\widetilde{e_\alpha})\widehat{c}(\widetilde{e_\beta})\Big)(x_{0}) \nonumber\\
=&\sum_{k<l <\alpha<\beta}T(\widetilde{e_k},\widetilde{e_l},\widetilde{e_\alpha},\widetilde{e_\beta})\sum_{a,b,c,d }u_{a}v_{b}w_{c}z_{d}
\delta_{k}^{a}{\rm{Tr}}\Big[c(\widetilde{e_b})
\widehat{c}(\widetilde{e_c})\widehat{c}(\widetilde{e_d})
c(\widetilde{e_l})\widehat{c}(\widetilde{e_\alpha})\widehat{c}(\widetilde{e_\beta})\Big](x_{0}) \nonumber\\
&-\sum_{k<l <\alpha<\beta}T(\widetilde{e_k},\widetilde{e_l},\widetilde{e_\alpha},\widetilde{e_\beta})\sum_{a,b,c,d }u_{a}v_{b}w_{c}z_{d}
\delta_{k}^{b}{\rm{Tr}}\Big[c(\widetilde{e_a})
\widehat{c}(\widetilde{e_c})\widehat{c}(\widetilde{e_d})
c(\widetilde{e_l})\widehat{c}(\widetilde{e_\alpha})\widehat{c}(\widetilde{e_\beta})\Big](x_{0}) \nonumber\\
=&-\frac{1}{6}T(u,v,w,z){\rm{Tr}}(\rm{Id}),
\end{align}
where we have used
 \begin{align}
 &\sum_{k<l <\alpha<\beta}T(\widetilde{e_k},\widetilde{e_l},\widetilde{e_\alpha},\widetilde{e_\beta})\sum_{a,b,c,d }u_{a}v_{b}w_{c}z_{d}
\delta_{k}^{a}{\rm{Tr}}\Big[c(\widetilde{e_b})
\widehat{c}(\widetilde{e_c})\widehat{c}(\widetilde{e_d})
c(\widetilde{e_l})\widehat{c}(\widetilde{e_\alpha})\widehat{c}(\widetilde{e_\beta})\Big](x_{0}) \nonumber\\
=&\sum_{k<l <\alpha<\beta}T(\widetilde{e_k},\widetilde{e_l},\widetilde{e_\alpha},\widetilde{e_\beta})\sum_{a,b,c,d }u_{a}v_{b}w_{c}z_{d}
\delta_{k}^{a}\delta_{l}^{b}{\rm{Tr}}\Big[
\widehat{c}(\widetilde{e_c})\widehat{c}(\widetilde{e_d})
 \widehat{c}(\widetilde{e_\alpha})\widehat{c}(\widetilde{e_\beta})\Big](x_{0}) \nonumber\\
 =&\sum_{k<l <\alpha<\beta}T(\widetilde{e_k},\widetilde{e_l},\widetilde{e_\alpha},\widetilde{e_\beta})\sum_{a,b,c,d }u_{a}v_{b}w_{c}z_{d}
\delta_{k}^{a}\delta_{l}^{b}\delta_{\alpha}^{d}\delta_{\beta}^{c}{\rm{Tr}}(\rm{Id}) \nonumber\\
 &-\sum_{k<l <\alpha<\beta}T(\widetilde{e_k},\widetilde{e_l},\widetilde{e_\alpha},\widetilde{e_\beta})\sum_{a,b,c,d }u_{a}v_{b}w_{c}z_{d}
\delta_{k}^{a}\delta_{l}^{b}\delta_{\alpha}^{c}\delta_{\beta}^{d}{\rm{Tr}}(\rm{Id}) \nonumber\\
 =&\frac{1}{4!}\Big(T(u,v,z,w)-T(u,v,w,z)\Big){\rm{Tr}}(\rm{Id}) \nonumber\\
 =&-\frac{1}{12}T(u,v,w,z){\rm{Tr}}(\rm{Id}),
\end{align}
and
 \begin{align}
 &\sum_{k<l <\alpha<\beta}T(\widetilde{e_k},\widetilde{e_l},\widetilde{e_\alpha},\widetilde{e_\beta})\sum_{a,b,c,d }u_{a}v_{b}w_{c}z_{d}
\delta_{k}^{b}{\rm{Tr}}\Big[c(\widetilde{e_a})
\widehat{c}(\widetilde{e_c})\widehat{c}(\widetilde{e_d})
c(\widetilde{e_l})\widehat{c}(\widetilde{e_\alpha})\widehat{c}(\widetilde{e_\beta})\Big](x_{0}) \nonumber\\
=&\sum_{k<l <\alpha<\beta}T(\widetilde{e_k},\widetilde{e_l},\widetilde{e_\alpha},\widetilde{e_\beta})\sum_{a,b,c,d }u_{a}v_{b}w_{c}z_{d}
\delta_{k}^{b}\delta_{l}^{a}{\rm{Tr}}\Big[
\widehat{c}(\widetilde{e_c})\widehat{c}(\widetilde{e_d})
 \widehat{c}(\widetilde{e_\alpha})\widehat{c}(\widetilde{e_\beta})\Big](x_{0}) \nonumber\\
 =&\sum_{k<l <\alpha<\beta}T(\widetilde{e_k},\widetilde{e_l},\widetilde{e_\alpha},\widetilde{e_\beta})\sum_{a,b,c,d }u_{a}v_{b}w_{c}z_{d}
\delta_{k}^{b}\delta_{l}^{a}\delta_{\alpha}^{d}\delta_{\beta}^{c}{\rm{Tr}}(\rm{Id}) \nonumber\\
 &-\sum_{k<l <\alpha<\beta}T(\widetilde{e_k},\widetilde{e_l},\widetilde{e_\alpha},\widetilde{e_\beta})\sum_{a,b,c,d }u_{a}v_{b}w_{c}z_{d}
\delta_{k}^{b}\delta_{l}^{a}\delta_{\alpha}^{c}\delta_{\beta}^{d}{\rm{Tr}}(\rm{Id}) \nonumber\\
 =&\frac{1}{4!}\Big(T(v,u,z,w)-T(v,u,w,z)\Big){\rm{Tr}}(\rm{Id}) \nonumber\\
 =&-\frac{1}{12}T(u,v,w,z){\rm{Tr}}(\rm{Id}).
\end{align}
Then the proof of the lemma is complete.
\end{proof}

 \begin{lem}
The following identities hold:
 \begin{align}
 &\int_{|\xi|=1}{\rm{Tr}}\Big(\sum_{i,j}g^{ij}
  \sum_{k<l <\alpha<\beta}T(\widetilde{e_k},\widetilde{e_l},\widetilde{e_\alpha},\widetilde{e_\beta})c(u)c(v)\widehat{c}(w)\widehat{c}(z)
  c(\partial_{i}) c(\widetilde{e_k})c(\widetilde{e_l})\widehat{c}(\widetilde{e_\alpha})\widehat{c}(\widetilde{e_\beta})
 \xi_{j} c(\xi) \Big)\sigma(\xi)(x_{0})\nonumber\\
 =&-\frac{1}{2}T(u,v,w,z){\rm{Tr}}(\rm{Id}){\rm{vol}}(S^{2m-1});\\
 &\int_{|\xi|=1}{\rm{Tr}}\Big(\sum_{i,j}g^{ij}
  \sum_{k<l <\alpha<\beta}T(\widetilde{e_k},\widetilde{e_l},\widetilde{e_\alpha},\widetilde{e_\beta})c(u)c(v)\widehat{c}(w)\widehat{c}(z)
 c(\widetilde{e_k})c(\widetilde{e_l})\widehat{c}(\widetilde{e_\alpha})\widehat{c}(\widetilde{e_\beta})
 c(\partial_{i})\xi_{j} c(\xi) \Big)\sigma(\xi)(x_{0})\nonumber\\
 =&\frac{1}{6}T(u,v,w,z){\rm{Tr}}(\rm{Id}){\rm{vol}}(S^{2m-1}).
\end{align}
\end{lem}
\begin{proof}
By the relation of the Clifford action and $ {\rm{Tr}}(AB)= {\rm{Tr}}(BA) $, in terms of the condition (2.24), we get
 \begin{align}
 &\int_{|\xi|=1}{\rm{Tr}}\Big(\sum_{i,j}g^{ij}
  \sum_{k<l <\alpha<\beta}T(\widetilde{e_k},\widetilde{e_l},\widetilde{e_\alpha},\widetilde{e_\beta})c(u)c(v)\widehat{c}(w)\widehat{c}(z)
  c(\partial_{i}) c(\widetilde{e_k})c(\widetilde{e_l})\widehat{c}(\widetilde{e_\alpha})\widehat{c}(\widetilde{e_\beta})
 \xi_{j} c(\xi) \Big)\sigma(\xi)(x_{0})\nonumber\\
 = &\int_{|\xi|=1}\sum_{i }\sum_{k<l <\alpha<\beta}T(\widetilde{e_k},\widetilde{e_l},\widetilde{e_\alpha},\widetilde{e_\beta})
 {\rm{Tr}}\Big( c(u)c(v)\widehat{c}(w)\widehat{c}(z)
 c(\widetilde{e_i})c(\widetilde{e_k})c(\widetilde{e_l})\widehat{c}(\widetilde{e_\alpha})\widehat{c}(\widetilde{e_\beta})
\xi_{i} \sum_{j }\xi_{j}c(\widetilde{e_j}) \Big)\sigma(\xi)(x_{0})\nonumber\\
  = &\sum_{i,j }\sum_{k<l <\alpha<\beta}\frac{1}{2m}\delta_{i}^{j}{\rm{vol}}(S^{2m-1})
  T(\widetilde{e_k},\widetilde{e_l},\widetilde{e_\alpha},\widetilde{e_\beta})
 {\rm{Tr}}\Big( c(u)c(v)\widehat{c}(w)\widehat{c}(z)
 c(\widetilde{e_i})c(\widetilde{e_k})c(\widetilde{e_l})\widehat{c}(\widetilde{e_\alpha})\widehat{c}(\widetilde{e_\beta})
c(\widetilde{e_j}) \Big)\sigma(\xi)(x_{0})\nonumber\\
  = & \sum_{i}\sum_{k<l <\alpha<\beta}\frac{1}{2m} {\rm{vol}}(S^{2m-1})
  T(\widetilde{e_k},\widetilde{e_l},\widetilde{e_\alpha},\widetilde{e_\beta})
 {\rm{Tr}}\Big( c(u)c(v)\widehat{c}(w)\widehat{c}(z)
 c(\widetilde{e_i})c(\widetilde{e_k})c(\widetilde{e_l})\widehat{c}(\widetilde{e_\alpha})\widehat{c}(\widetilde{e_\beta})
c(\widetilde{e_i}) \Big)\sigma(\xi)(x_{0})\nonumber\\
  = & \sum_{i}\sum_{k<l <\alpha<\beta}  \sum_{a,b,c,d }\frac{1}{2m} {\rm{vol}}(S^{2m-1})
  T(\widetilde{e_k},\widetilde{e_l},\widetilde{e_\alpha},\widetilde{e_\beta})u_{a}v_{b}w_{c}z_{d}\nonumber\\
 &\times{\rm{Tr}}\Big( c(\widetilde{e_a})c(\widetilde{e_b})\widehat{c}(\widetilde{e_c})\widehat{c}(\widetilde{e_d})
 c(\widetilde{e_i})c(\widetilde{e_k})c(\widetilde{e_l})\widehat{c}(\widetilde{e_\alpha})\widehat{c}(\widetilde{e_\beta})
c(\widetilde{e_i}) \Big)\sigma(\xi)(x_{0})\nonumber\\
  = & \sum_{i}\sum_{k<l <\alpha<\beta}  \sum_{a,b,c,d }\frac{1}{2m} {\rm{vol}}(S^{2m-1})
  T(\widetilde{e_k},\widetilde{e_l},\widetilde{e_\alpha},\widetilde{e_\beta})u_{a}v_{b}w_{c}z_{d}\nonumber\\
 &\times{\rm{Tr}}\Big[ c(\widetilde{e_a})\Big(-c(\widetilde{e_i})c(\widetilde{e_b})-2\delta_{b}^{i}\Big)
 \widehat{c}(\widetilde{e_c})\widehat{c}(\widetilde{e_d})
c(\widetilde{e_k})c(\widetilde{e_l})\widehat{c}(\widetilde{e_\alpha})\widehat{c}(\widetilde{e_\beta})
c(\widetilde{e_i}) \Big]\sigma(\xi)(x_{0})\nonumber\\
  = &\cdots\nonumber\\
  = & 3\sum_{k<l <\alpha<\beta}{\rm{vol}}(S^{2m-1})
  T(\widetilde{e_k},\widetilde{e_l},\widetilde{e_\alpha},\widetilde{e_\beta})
 {\rm{Tr}}\Big( c(u)c(v)\widehat{c}(w)\widehat{c}(z)
 c(\widetilde{e_k})c(\widetilde{e_l})\widehat{c}(\widetilde{e_\alpha})\widehat{c}(\widetilde{e_\beta})
  \Big)\sigma(\xi)(x_{0})\nonumber\\
    = &-\frac{1}{2}T(u,v,w,z){\rm{Tr}}(\rm{Id}).
\end{align}
On the other hand,
 \begin{align}
 &\int_{|\xi|=1}{\rm{Tr}}\Big(\sum_{i,j}g^{ij}
  \sum_{k<l <\alpha<\beta}T(\widetilde{e_k},\widetilde{e_l},\widetilde{e_\alpha},\widetilde{e_\beta})c(u)c(v)\widehat{c}(w)\widehat{c}(z)
 c(\widetilde{e_k})c(\widetilde{e_l})\widehat{c}(\widetilde{e_\alpha})\widehat{c}(\widetilde{e_\beta})
 c(\partial_{i})\xi_{j} c(\xi) \Big)\sigma(\xi)(x_{0})\nonumber\\
 = &\int_{|\xi|=1}\sum_{i }\sum_{k<l <\alpha<\beta}T(\widetilde{e_k},\widetilde{e_l},\widetilde{e_\alpha},\widetilde{e_\beta})
 {\rm{Tr}}\Big( c(u)c(v)\widehat{c}(w)\widehat{c}(z)
 c(\widetilde{e_k})c(\widetilde{e_l})\widehat{c}(\widetilde{e_\alpha})\widehat{c}(\widetilde{e_\beta})
 c(\widetilde{e_i})\xi_{i} \sum_{j }\xi_{j}c(\widetilde{e_j}) \Big)\sigma(\xi)(x_{0})\nonumber\\
  = &\sum_{i,j }\sum_{k<l <\alpha<\beta}\frac{1}{2m}\delta_{i}^{j}{\rm{vol}}(S^{2m-1})
  T(\widetilde{e_k},\widetilde{e_l},\widetilde{e_\alpha},\widetilde{e_\beta})
 {\rm{Tr}}\Big( c(u)c(v)\widehat{c}(w)\widehat{c}(z)
 c(\widetilde{e_k})c(\widetilde{e_l})\widehat{c}(\widetilde{e_\alpha})\widehat{c}(\widetilde{e_\beta})
 c(\widetilde{e_i})c(\widetilde{e_j}) \Big)\sigma(\xi)(x_{0})\nonumber\\
   = &-\sum_{i}\sum_{k<l <\alpha<\beta}\frac{1}{2m}{\rm{vol}}(S^{2m-1})
  T(\widetilde{e_k},\widetilde{e_l},\widetilde{e_\alpha},\widetilde{e_\beta})
 {\rm{Tr}}\Big( c(u)c(v)\widehat{c}(w)\widehat{c}(z)
 c(\widetilde{e_k})c(\widetilde{e_l})\widehat{c}(\widetilde{e_\alpha})\widehat{c}(\widetilde{e_\beta})
  \Big)\sigma(\xi)(x_{0})\nonumber\\
  =&-\frac{1}{6}T(u,v,w,z){\rm{vol}}(S^{2m-1}){\rm{Tr}}(\rm{Id}).
\end{align}
Then the proof of the lemma is complete.
\end{proof}
From (4.3), Lemma 4.5 and Lemma 4.6, then  the proof of Theorem 1.4 is complete.

\subsection{The proof of Theorem 1.5}
From (2.18) and (4.2), we obtain
\begin{align}
 &\sigma_{-2m}(\widehat{c}(u)\widehat{c}(v)\widehat{c}(w)\widehat{c}(z)(d+\delta)_{T_{4}}^{1-2m})(x_{0})|_{|\xi|=1} \nonumber\\
 =&\widehat{c}(u)\widehat{c}(v)\widehat{c}(w)\widehat{c}(z)\Big\{|\xi|^{-2m}\sigma_0((d+\delta)_{T_{4}})
 +\sum_{j=1}^{2m+2}\partial_{\xi_j}(|\xi|^{-2m})\partial_{x_j}c(\xi)+
\Big[m\sigma_2((d+\delta)_{T_{4}})^{-m+1}\sigma_{-3}((d+\delta)_{T_{4}}^{-2})  \nonumber\\
&-\sqrt{-1} \sum_{k=0}^{m-2}\sum_{\mu=1}^{2m+2}
\partial_{\xi_{\mu}}\sigma_{2}^{-m+k+1}((d+\delta)_{T_{4}}^2)
\partial_{x_{\mu}}\sigma_{2}^{-1}((d+\delta)_{T_{4}}^2)(\sigma_2((d+\delta)_{T_{4}}^2))^{-k}\Big ]\sqrt{-1}c(\xi)\Big\}(x_{0})|_{|\xi|=1} \nonumber\\
=&\widehat{c}(u)\widehat{c}(v)\widehat{c}(w)\widehat{c}(z)\sqrt{-1}\sum_{k<l <\alpha<\beta}
T(\widetilde{e_k},\widetilde{e_l},\widetilde{e_\alpha},\widetilde{e_\beta})
\widehat{c}(\widetilde{e_k})\widehat{c}(\widetilde{e_l})\widehat{c}(\widetilde{e_\alpha})\widehat{c}(\widetilde{e_\beta}) \nonumber\\
&+m\sum_{i,j}g^{ij}\widehat{c}(u)\widehat{c}(v)\widehat{c}(w)\widehat{c}(z)\sqrt{-1}
\Big[c(\partial_{i}) \sum_{k<l <\alpha<\beta}T(\widetilde{e_k},\widetilde{e_l},\widetilde{e_\alpha},\widetilde{e_\beta})
\widehat{c}(\widetilde{e_k})\widehat{c}(\widetilde{e_l})\widehat{c}(\widetilde{e_\alpha})\widehat{c}(\widetilde{e_\beta})
 \nonumber\\
&+ \sum_{k<l <\alpha<\beta}T(\widetilde{e_k},\widetilde{e_l},\widetilde{e_\alpha},\widetilde{e_\beta})
\widehat{c}(\widetilde{e_k})\widehat{c}(\widetilde{e_l})\widehat{c}(\widetilde{e_\alpha})\widehat{c}(\widetilde{e_\beta}) c(\partial_{i})   \Big]\xi_{j}c(\xi).
\end{align}

Let $u=\sum_{i=1}^{n}u_{i}\widetilde{e_{i}}$, $v=\sum_{j=1}^{n}v_{j}\widetilde{e_{j}} $, $w=\sum_{k=1}^{n}w_{k}\widetilde{e_{k}} $,
$z=\sum_{l=1}^{n}z_{l}\widetilde{e_{l}} $,
where $\{\widetilde{e_{1}},\widetilde{e_{2}},\cdots,\widetilde{e_{n}}\}$ is the orthogonal basis about $g^{TM}$,
then  $\widehat{c}(u)= \sum_{i=1}^{n}u_{i}\widehat{c}(\widetilde{e_{i}})$, $\widehat{c}(v)= \sum_{j=1}^{n}v_{j}\widehat{c}(\widetilde{e_{j}})$,
$\widehat{c}(w)= \sum_{k=1}^{n}w_{k}\widehat{c}(\widetilde{e_{k}})$, $\widehat{c}(z)= \sum_{l=1}^{n}z_{l}\widehat{c}(\widetilde{e_{l}})$.
The following Lemmas of traces in terms of the Clifford action
is very efficient for solving the spectral torsion  associated with the de-Rham Hodge type operator.
 \begin{lem}
The following identities hold:
 \begin{align}
 {\rm{Tr}}\Big(\sum_{k<l <\alpha<\beta}T(\widetilde{e_k},\widetilde{e_l},\widetilde{e_\alpha},\widetilde{e_\beta})
 \widehat{c}(u)\widehat{c}(v)\widehat{c}(w)\widehat{c}(z)
\widehat{c}(\widetilde{e_k})\widehat{c}(\widetilde{e_l})\widehat{c}(\widetilde{e_\alpha})\widehat{c}(\widetilde{e_\beta})\Big)(x_{0})
= T(u,v,w,z){\rm{Tr}}(\rm{Id}).
\end{align}
\end{lem}
\begin{proof}
By the relation of the Clifford action and $ {\rm{Tr}}(AB)= {\rm{Tr}}(BA) $, we get
 \begin{align}
 &{\rm{Tr}}\Big(\sum_{k<l <\alpha<\beta}T(\widetilde{e_k},\widetilde{e_l},\widetilde{e_\alpha},\widetilde{e_\beta})
\widehat{c}(u)\widehat{c}(v)\widehat{c}(w)\widehat{c}(z)
\widehat{c}(\widetilde{e_k})\widehat{c}(\widetilde{e_l})\widehat{c}(\widetilde{e_\alpha})\widehat{c}(\widetilde{e_\beta})\Big)(x_{0}) \nonumber\\
= &\sum_{k<l <\alpha<\beta}T(\widetilde{e_k},\widetilde{e_l},\widetilde{e_\alpha},\widetilde{e_\beta})\sum_{a,b,c,d }u_{a}v_{b}w_{c}z_{d}
{\rm{Tr}}\Big(\widehat{c}(\widetilde{e_a})\widehat{c}(\widetilde{e_b})\widehat{c}(\widetilde{e_c})\widehat{c}(\widetilde{e_d})
\widehat{c}(\widetilde{e_k})\widehat{c}(\widetilde{e_l})\widehat{c}(\widetilde{e_\alpha})\widehat{c}(\widetilde{e_\beta})\Big)(x_{0}) \nonumber\\
= &\sum_{k<l <\alpha<\beta}T(\widetilde{e_k},\widetilde{e_l},\widetilde{e_\alpha},\widetilde{e_\beta})\sum_{a,b,c,d }u_{a}v_{b}w_{c}z_{d}
{\rm{Tr}}\Big[\widehat{c}(\widetilde{e_a})\widehat{c}(\widetilde{e_b})\widehat{c}(\widetilde{e_c})
\Big(-\widehat{c}(\widetilde{e_k})\widehat{c}(\widetilde{e_d})+2\delta_{k}^{d}   \Big)
\widehat{c}(\widetilde{e_l})\widehat{c}(\widetilde{e_\alpha})\widehat{c}(\widetilde{e_\beta})\Big](x_{0}) \nonumber\\
= &-\sum_{k<l <\alpha<\beta}T(\widetilde{e_k},\widetilde{e_l},\widetilde{e_\alpha},\widetilde{e_\beta})\sum_{a,b,c,d }u_{a}v_{b}w_{c}z_{d}
{\rm{Tr}}\Big[\widehat{c}(\widetilde{e_a}) \widehat{c}(\widetilde{e_b})\widehat{c}(\widetilde{e_c})
\widehat{c}(\widetilde{e_k})\widehat{c}(\widetilde{e_d})
\widehat{c}(\widetilde{e_l})\widehat{c}(\widetilde{e_\alpha})\widehat{c}(\widetilde{e_\beta})\Big](x_{0}) \nonumber\\
&+\sum_{k<l <\alpha<\beta}T(\widetilde{e_k},\widetilde{e_l},\widetilde{e_\alpha},\widetilde{e_\beta})\sum_{a,b,c,d }u_{a}v_{b}w_{c}z_{d}
2\delta_{k}^{d}{\rm{Tr}}\Big[\widehat{c}(\widetilde{e_a})
\widehat{c}(\widetilde{e_b})\widehat{c}(\widetilde{e_c})
\widehat{c}(\widetilde{e_l})\widehat{c}(\widetilde{e_\alpha})\widehat{c}(\widetilde{e_\beta})\Big](x_{0}) \nonumber\\
= &\cdots\nonumber\\
= &-\sum_{k<l <\alpha<\beta}T(\widetilde{e_k},\widetilde{e_l},\widetilde{e_\alpha},\widetilde{e_\beta})\sum_{a,b,c,d }u_{a}v_{b}w_{c}z_{d}
{\rm{Tr}}\Big[\widehat{c}(\widetilde{e_a})\widehat{c}(\widetilde{e_b})
\widehat{c}(\widetilde{e_c})\widehat{c}(\widetilde{e_d})
\widehat{c}(\widetilde{e_k})\widehat{c}(\widetilde{e_l})  \widehat{c}(\widetilde{e_\alpha})\widehat{c}(\widetilde{e_\beta})\Big](x_{0}) \nonumber\\
&-\sum_{k<l <\alpha<\beta}T(\widetilde{e_k},\widetilde{e_l},\widetilde{e_\alpha},\widetilde{e_\beta})\sum_{a,b,c,d }u_{a}v_{b}w_{c}z_{d}
2\delta_{k}^{a}{\rm{Tr}}\Big[\widehat{c}(\widetilde{e_b})
\widehat{c}(\widetilde{e_c})\widehat{c}(\widetilde{e_d})
\widehat{c}(\widetilde{e_l})\widehat{c}(\widetilde{e_\alpha})\widehat{c}(\widetilde{e_\beta})\Big](x_{0}) \nonumber\\
&+\sum_{k<l <\alpha<\beta}T(\widetilde{e_k},\widetilde{e_l},\widetilde{e_\alpha},\widetilde{e_\beta})\sum_{a,b,c,d }u_{a}v_{b}w_{c}z_{d}
2\delta_{k}^{b}{\rm{Tr}}\Big[\widehat{c}(\widetilde{e_a})
\widehat{c}(\widetilde{e_c})\widehat{c}(\widetilde{e_d})
\widehat{c}(\widetilde{e_l})\widehat{c}(\widetilde{e_\alpha})\widehat{c}(\widetilde{e_\beta})\Big](x_{0}) \nonumber\\
&-\sum_{k<l <\alpha<\beta}T(\widetilde{e_k},\widetilde{e_l},\widetilde{e_\alpha},\widetilde{e_\beta})\sum_{a,b,c,d }u_{a}v_{b}w_{c}z_{d}
2\delta_{k}^{c}{\rm{Tr}}\Big[\widehat{c}(\widetilde{e_a})
\widehat{c}(\widetilde{e_b})\widehat{c}(\widetilde{e_d})
\widehat{c}(\widetilde{e_l})\widehat{c}(\widetilde{e_\alpha})\widehat{c}(\widetilde{e_\beta})\Big](x_{0}) \nonumber\\
&+\sum_{k<l <\alpha<\beta}T(\widetilde{e_k},\widetilde{e_l},\widetilde{e_\alpha},\widetilde{e_\beta})\sum_{a,b,c,d }u_{a}v_{b}w_{c}z_{d}
2\delta_{k}^{d}{\rm{Tr}}\Big[\widehat{c}(\widetilde{e_a})
\widehat{c}(\widetilde{e_b})\widehat{c}(\widetilde{e_c})
\widehat{c}(\widetilde{e_l})\widehat{c}(\widetilde{e_\alpha})\widehat{c}(\widetilde{e_\beta})\Big](x_{0}),
\end{align}

where we have used
 \begin{align}
 &\sum_{k<l <\alpha<\beta}T(\widetilde{e_k},\widetilde{e_l},\widetilde{e_\alpha},\widetilde{e_\beta})\sum_{a,b,c,d }u_{a}v_{b}w_{c}z_{d}
\delta_{k}^{a}{\rm{Tr}}\Big[c(\widetilde{e_b})
\widehat{c}(\widetilde{e_c})\widehat{c}(\widetilde{e_d})
\widehat{c}(\widetilde{e_l})\widehat{c}(\widetilde{e_\alpha})\widehat{c}(\widetilde{e_\beta})\Big](x_{0}) \nonumber\\
=&\sum_{k<l <\alpha<\beta}T(\widetilde{e_k},\widetilde{e_l},\widetilde{e_\alpha},\widetilde{e_\beta})\sum_{a,b,c,d }u_{a}v_{b}w_{c}z_{d}
\delta_{k}^{a}{\rm{Tr}}\Big[c(\widetilde{e_b})
\widehat{c}(\widetilde{e_c})\Big(-\widehat{c}(\widetilde{e_l})\widehat{c}(\widetilde{e_d})+2\delta_{l}^{d} \Big)
\widehat{c}(\widetilde{e_\alpha})\widehat{c}(\widetilde{e_\beta})\Big](x_{0}) \nonumber\\
=&-\sum_{k<l <\alpha<\beta}T(\widetilde{e_k},\widetilde{e_l},\widetilde{e_\alpha},\widetilde{e_\beta})\sum_{a,b,c,d }u_{a}v_{b}w_{c}z_{d}
\delta_{k}^{a}{\rm{Tr}}\Big[c(\widetilde{e_b})
\widehat{c}(\widetilde{e_c}) \widehat{c}(\widetilde{e_l})\widehat{c}(\widetilde{e_d})
\widehat{c}(\widetilde{e_\alpha})\widehat{c}(\widetilde{e_\beta})\Big](x_{0}) \nonumber\\
&+\sum_{k<l <\alpha<\beta}T(\widetilde{e_k},\widetilde{e_l},\widetilde{e_\alpha},\widetilde{e_\beta})\sum_{a,b,c,d }u_{a}v_{b}w_{c}z_{d}
2\delta_{k}^{a}\delta_{l}^{d}{\rm{Tr}}\Big[c(\widetilde{e_b})
\widehat{c}(\widetilde{e_c})
\widehat{c}(\widetilde{e_\alpha})\widehat{c}(\widetilde{e_\beta})\Big](x_{0}) \nonumber\\
= &\cdots\nonumber\\
=&-\sum_{k<l <\alpha<\beta}T(\widetilde{e_k},\widetilde{e_l},\widetilde{e_\alpha},\widetilde{e_\beta})\sum_{a,b,c,d }u_{a}v_{b}w_{c}z_{d}
\delta_{k}^{a}{\rm{Tr}}\Big[c(\widetilde{e_b})
\widehat{c}(\widetilde{e_c})\widehat{c}(\widetilde{e_d})
\widehat{c}(\widetilde{e_l})\widehat{c}(\widetilde{e_\alpha})\widehat{c}(\widetilde{e_\beta})\Big](x_{0}) \nonumber\\
&+\sum_{k<l <\alpha<\beta}T(\widetilde{e_k},\widetilde{e_l},\widetilde{e_\alpha},\widetilde{e_\beta})\sum_{a,b,c,d }u_{a}v_{b}w_{c}z_{d}
2\delta_{k}^{a}\delta_{l}^{d}{\rm{Tr}}\Big[c(\widetilde{e_b})
\widehat{c}(\widetilde{e_c})
\widehat{c}(\widetilde{e_\alpha})\widehat{c}(\widetilde{e_\beta})\Big](x_{0}) \nonumber\\
&-\sum_{k<l <\alpha<\beta}T(\widetilde{e_k},\widetilde{e_l},\widetilde{e_\alpha},\widetilde{e_\beta})\sum_{a,b,c,d }u_{a}v_{b}w_{c}z_{d}
2\delta_{k}^{a}\delta_{l}^{c}{\rm{Tr}}\Big[c(\widetilde{e_b})
\widehat{c}(\widetilde{e_d})
\widehat{c}(\widetilde{e_\alpha})\widehat{c}(\widetilde{e_\beta})\Big](x_{0}) \nonumber\\
&+\sum_{k<l <\alpha<\beta}T(\widetilde{e_k},\widetilde{e_l},\widetilde{e_\alpha},\widetilde{e_\beta})\sum_{a,b,c,d }u_{a}v_{b}w_{c}z_{d}
2\delta_{k}^{a}\delta_{l}^{b}{\rm{Tr}}\Big[c(\widetilde{e_c})
\widehat{c}(\widetilde{e_d})
\widehat{c}(\widetilde{e_\alpha})\widehat{c}(\widetilde{e_\beta})\Big](x_{0}).
\end{align}

By the generating relation, we see that the left hand side of (4.18) equals
 \begin{align}
 &\sum_{k<l <\alpha<\beta}T(\widetilde{e_k},\widetilde{e_l},\widetilde{e_\alpha},\widetilde{e_\beta})\sum_{a,b,c,d }u_{a}v_{b}w_{c}z_{d}
\delta_{k}^{a}{\rm{Tr}}\Big[c(\widetilde{e_b})
\widehat{c}(\widetilde{e_c})\widehat{c}(\widetilde{e_d})
\widehat{c}(\widetilde{e_l})\widehat{c}(\widetilde{e_\alpha})\widehat{c}(\widetilde{e_\beta})\Big](x_{0}) \nonumber\\
=&\sum_{k<l <\alpha<\beta}T(\widetilde{e_k},\widetilde{e_l},\widetilde{e_\alpha},\widetilde{e_\beta})\sum_{a,b,c,d }u_{a}v_{b}w_{c}z_{d}
\delta_{k}^{a}\delta_{l}^{d}{\rm{Tr}}\Big[c(\widetilde{e_b})
\widehat{c}(\widetilde{e_c})
\widehat{c}(\widetilde{e_\alpha})\widehat{c}(\widetilde{e_\beta})\Big](x_{0}) \nonumber\\
&-\sum_{k<l <\alpha<\beta}T(\widetilde{e_k},\widetilde{e_l},\widetilde{e_\alpha},\widetilde{e_\beta})\sum_{a,b,c,d }u_{a}v_{b}w_{c}z_{d}
\delta_{k}^{a}\delta_{l}^{c}{\rm{Tr}}\Big[c(\widetilde{e_b})
\widehat{c}(\widetilde{e_d})
\widehat{c}(\widetilde{e_\alpha})\widehat{c}(\widetilde{e_\beta})\Big](x_{0}) \nonumber\\
&+\sum_{k<l <\alpha<\beta}T(\widetilde{e_k},\widetilde{e_l},\widetilde{e_\alpha},\widetilde{e_\beta})\sum_{a,b,c,d }u_{a}v_{b}w_{c}z_{d}
\delta_{k}^{a}\delta_{l}^{b}{\rm{Tr}}\Big[c(\widetilde{e_c})
\widehat{c}(\widetilde{e_d})
\widehat{c}(\widetilde{e_\alpha})\widehat{c}(\widetilde{e_\beta})\Big](x_{0})   \nonumber\\
=&\frac{1}{24}\big(T(u,v,z,w)-T(u,v,w,z)\big){\rm{Tr}}(\rm{Id})
-\frac{1}{24}\big(T(u,w,z,v)-T(u,w,v,z)\big){\rm{Tr}}(\rm{Id})\nonumber\\
&+\frac{1}{24}\big(T(u,z,w,v)-T(u,z,v,w)\big){\rm{Tr}}(\rm{Id})\nonumber\\
=&-\frac{1}{4}T(u,v,z,w){\rm{Tr}}(\rm{Id}).
\end{align}
Similarly

 \begin{align}
&\sum_{k<l <\alpha<\beta}T(\widetilde{e_k},\widetilde{e_l},\widetilde{e_\alpha},\widetilde{e_\beta})\sum_{a,b,c,d }u_{a}v_{b}w_{c}z_{d}
\delta_{k}^{b}{\rm{Tr}}\Big[\widehat{c}(\widetilde{e_a})
\widehat{c}(\widetilde{e_c})\widehat{c}(\widetilde{e_d})
\widehat{c}(\widetilde{e_l})\widehat{c}(\widetilde{e_\alpha})\widehat{c}(\widetilde{e_\beta})\Big](x_{0}) \nonumber\\
=&\frac{1}{4}T(u,v,z,w){\rm{Tr}}(\rm{Id}).
\end{align}

 \begin{align}
&\sum_{k<l <\alpha<\beta}T(\widetilde{e_k},\widetilde{e_l},\widetilde{e_\alpha},\widetilde{e_\beta})\sum_{a,b,c,d }u_{a}v_{b}w_{c}z_{d}
 \delta_{k}^{c}{\rm{Tr}}\Big[\widehat{c}(\widetilde{e_a})
\widehat{c}(\widetilde{e_b})\widehat{c}(\widetilde{e_d})
\widehat{c}(\widetilde{e_l})\widehat{c}(\widetilde{e_\alpha})\widehat{c}(\widetilde{e_\beta})\Big](x_{0}) \nonumber\\
=&-\frac{1}{4}T(u,v,z,w){\rm{Tr}}(\rm{Id}).
\end{align}

 \begin{align}
&\sum_{k<l <\alpha<\beta}T(\widetilde{e_k},\widetilde{e_l},\widetilde{e_\alpha},\widetilde{e_\beta})\sum_{a,b,c,d }u_{a}v_{b}w_{c}z_{d}
 \delta_{k}^{d}{\rm{Tr}}\Big[\widehat{c}(\widetilde{e_a})
\widehat{c}(\widetilde{e_b})\widehat{c}(\widetilde{e_c})
\widehat{c}(\widetilde{e_l})\widehat{c}(\widetilde{e_\alpha})\widehat{c}(\widetilde{e_\beta})\Big](x_{0}) \nonumber\\
=&\frac{1}{4}T(u,v,z,w){\rm{Tr}}(\rm{Id}).
\end{align}
Therefor, we obtain
\begin{align}
 &{\rm{Tr}}\Big(\sum_{k<l <\alpha<\beta}T(\widetilde{e_k},\widetilde{e_l},\widetilde{e_\alpha},\widetilde{e_\beta})
\widehat{c}(u)\widehat{c}(v)\widehat{c}(w)\widehat{c}(z)
\widehat{c}(\widetilde{e_k})\widehat{c}(\widetilde{e_l})\widehat{c}(\widetilde{e_\alpha})\widehat{c}(\widetilde{e_\beta})\Big)(x_{0}) \nonumber\\
=&T(u,v,z,w){\rm{Tr}}(\rm{Id}).
\end{align}
Then the proof of the lemma is complete.
\end{proof}

 \begin{lem}
The following identities hold:
 \begin{align}
 &\int_{|\xi|=1}{\rm{Tr}}\Big(\sum_{i,j}g^{ij}
  \sum_{k<l <\alpha<\beta}T(\widetilde{e_k},\widetilde{e_l},\widetilde{e_\alpha},\widetilde{e_\beta})
  \widehat{c}(u)\widehat{c}(v)\widehat{c}(w)\widehat{c}(z)
c(\partial_{i})\widehat{c}(\widetilde{e_k})c(\widetilde{e_l})\widehat{c}(\widetilde{e_\alpha})\widehat{c}(\widetilde{e_\beta})
 \xi_{j} c(\xi) \Big)\sigma(\xi)(x_{0})\nonumber\\
 =&-T(u,v,w,z){\rm{Tr}}(\rm{Id}){\rm{vol}}(S^{2m-1});\\
 &\int_{|\xi|=1}{\rm{Tr}}\Big(\sum_{i,j}g^{ij}
  \sum_{k<l <\alpha<\beta}T(\widetilde{e_k},\widetilde{e_l},\widetilde{e_\alpha},\widetilde{e_\beta})
  \widehat{c}(u)\widehat{c}(v)\widehat{c}(w)\widehat{c}(z)
 \widehat{c}(\widetilde{e_k})\widehat{c}(\widetilde{e_l})\widehat{c}(\widetilde{e_\alpha})\widehat{c}(\widetilde{e_\beta})
 c(\partial_{i})\xi_{j} c(\xi) \Big)\sigma(\xi)(x_{0})\nonumber\\
 =&-T(u,v,w,z){\rm{Tr}}(\rm{Id}){\rm{vol}}(S^{2m-1}).
\end{align}
\end{lem}
\begin{proof}
By the relation of the Clifford action and $ {\rm{Tr}}(AB)= {\rm{Tr}}(BA) $, in terms of the condition (2.24), we get
 \begin{align}
 &\int_{|\xi|=1}{\rm{Tr}}\Big(\sum_{i,j}g^{ij}
  \sum_{k<l <\alpha<\beta}T(\widetilde{e_k},\widetilde{e_l},\widetilde{e_\alpha},\widetilde{e_\beta})
  \widehat{c}(u)\widehat{c}(v)\widehat{c}(w)\widehat{c}(z)
  c(\partial_{i}) \widehat{c}(\widetilde{e_k})\widehat{c}(\widetilde{e_l})\widehat{c}(\widetilde{e_\alpha})\widehat{c}(\widetilde{e_\beta})
 \xi_{j} c(\xi) \Big)\sigma(\xi)(x_{0})\nonumber\\
 = &\int_{|\xi|=1}\sum_{i }\sum_{k<l <\alpha<\beta}T(\widetilde{e_k},\widetilde{e_l},\widetilde{e_\alpha},\widetilde{e_\beta})
 {\rm{Tr}}\Big( \widehat{c}(u)\widehat{c}(v)\widehat{c}(w)\widehat{c}(z)
 c(\widetilde{e_i})\widehat{c}(\widetilde{e_k})\widehat{c}(\widetilde{e_l})\widehat{c}(\widetilde{e_\alpha})\widehat{c}(\widetilde{e_\beta})
\xi_{i} \sum_{j }\xi_{j}c(\widetilde{e_j}) \Big)\sigma(\xi)(x_{0})\nonumber\\
  = &\sum_{i,j }\sum_{k<l <\alpha<\beta}\frac{1}{2m}\delta_{i}^{j}{\rm{vol}}(S^{2m-1})
  T(\widetilde{e_k},\widetilde{e_l},\widetilde{e_\alpha},\widetilde{e_\beta})
 {\rm{Tr}}\Big( \widehat{c}(u)\widehat{c}(v)\widehat{c}(w)\widehat{c}(z)
 c(\widetilde{e_i})\widehat{c}(\widetilde{e_k})\widehat{c}(\widetilde{e_l})\widehat{c}(\widetilde{e_\alpha})\widehat{c}(\widetilde{e_\beta})
c(\widetilde{e_j}) \Big)\sigma(\xi)(x_{0})\nonumber\\
    = & -\sum_{k<l <\alpha<\beta}{\rm{vol}}(S^{2m-1})
  T(\widetilde{e_k},\widetilde{e_l},\widetilde{e_\alpha},\widetilde{e_\beta})
 {\rm{Tr}}\Big( \widehat{c}(u)\widehat{c}(v)\widehat{c}(w)\widehat{c}(z)
 \widehat{c}(\widetilde{e_k})\widehat{c}(\widetilde{e_l})\widehat{c}(\widetilde{e_\alpha})\widehat{c}(\widetilde{e_\beta})
  \Big)\sigma(\xi)(x_{0})\nonumber\\
    = &-T(u,v,w,z){\rm{Tr}}(\rm{Id}).
\end{align}
On the other hand,
 \begin{align}
 &\int_{|\xi|=1}{\rm{Tr}}\Big(\sum_{i,j}g^{ij}
  \sum_{k<l <\alpha<\beta}T(\widetilde{e_k},\widetilde{e_l},\widetilde{e_\alpha},\widetilde{e_\beta})
  \widehat{c}(u)\widehat{c}(v)\widehat{c}(w)\widehat{c}(z)
  \widehat{c}(\widetilde{e_k})\widehat{c}(\widetilde{e_l})\widehat{c}(\widetilde{e_\alpha})\widehat{c}(\widetilde{e_\beta})
 c(\partial_{i})\xi_{j} c(\xi) \Big)\sigma(\xi)(x_{0})\nonumber\\
 = &\int_{|\xi|=1}\sum_{i }\sum_{k<l <\alpha<\beta}T(\widetilde{e_k},\widetilde{e_l},\widetilde{e_\alpha},\widetilde{e_\beta})
 {\rm{Tr}}\Big( \widehat{c}(u)\widehat{c}(v)\widehat{c}(w)\widehat{c}(z)
\widehat{c}(\widetilde{e_k})\widehat{c}(\widetilde{e_l})\widehat{c}(\widetilde{e_\alpha})\widehat{c}(\widetilde{e_\beta})
c(\widetilde{e_i})\xi_{i} \sum_{j }\xi_{j}c(\widetilde{e_j}) \Big)\sigma(\xi)(x_{0})\nonumber\\
  = &\sum_{i,j }\sum_{k<l <\alpha<\beta}\frac{1}{2m}\delta_{i}^{j}{\rm{vol}}(S^{2m-1})
  T(\widetilde{e_k},\widetilde{e_l},\widetilde{e_\alpha},\widetilde{e_\beta})
 {\rm{Tr}}\Big( \widehat{c}(u)\widehat{c}(v)\widehat{c}(w)\widehat{c}(z)
\widehat{c}(\widetilde{e_k})\widehat{c}(\widetilde{e_l})\widehat{c}(\widetilde{e_\alpha})\widehat{c}(\widetilde{e_\beta})
 c(\widetilde{e_i})c(\widetilde{e_j}) \Big)\sigma(\xi)(x_{0})\nonumber\\
    = & -\sum_{k<l <\alpha<\beta}{\rm{vol}}(S^{2m-1})
  T(\widetilde{e_k},\widetilde{e_l},\widetilde{e_\alpha},\widetilde{e_\beta})
 {\rm{Tr}}\Big( \widehat{c}(u)\widehat{c}(v)\widehat{c}(w)\widehat{c}(z)
 \widehat{c}(\widetilde{e_k})\widehat{c}(\widetilde{e_l})\widehat{c}(\widetilde{e_\alpha})\widehat{c}(\widetilde{e_\beta})
  \Big)\sigma(\xi)(x_{0})\nonumber\\
    = &-T(u,v,w,z){\rm{Tr}}(\rm{Id}).
\end{align}
Then the proof of the lemma is complete.
\end{proof}
From (4.5), Lemma 4.7 and lemma 4.8, then  the proof of Theorem 1.5 is complete.

\section{The spectral torsion for de-Rham Hodge operator on manifold with boundary }
The purpose of this section is to specify the spectral torsion on manifold with boundary  for the de-Rham Hodge type operator with torsion
$(d+\delta)_{T_{2}}$. Now we recall the main theorem in \cite{FGLS}.
\begin{thm}\cite{FGLS}
 Let $X$ and $\partial X$ be connected, ${\rm dim}X=n\geq3$,
 $A=\left(\begin{array}{lcr}\pi^+P+G &   K \\
T &  S    \end{array}\right)$ $\in \mathcal{B}$ , and denote by $p$, $b$ and $s$ the local symbols of $P,G$ and $S$ respectively.
 Define:
 \begin{align}
{\rm{\widetilde{Wres}}}(A)=&\int_X\int_{\bf S}{\mathrm{Tr}}_E\left[p_{-n}(x,\xi)\right]\sigma(\xi){\rm d}x \nonumber\\
&+2\pi\int_ {\partial X}\int_{\bf S'}\left\{{\mathrm{Tr}}_E\left[({\mathrm{Tr}}b_{-n})(x',\xi')\right]+{\mathrm{Tr}}
_F\left[s_{1-n}(x',\xi')\right]\right\}\sigma(\xi'){\rm d}x',
\end{align}
Then~~ a) ${\rm \widetilde{Wres}}([A,B])=0 $, for any
$A,B\in\mathcal{B}$;~~ b) It is a unique continuous trace on
$\mathcal{B}/\mathcal{B}^{-\infty}$.
\end{thm}
 Let $M$ be a compact oriented Riemannian manifold of even dimension $n=2m$. Let $p_{1},p_{2}$ be nonnegative integers and $p_{1}+p_{2}\leq n$.
 Denote by $\sigma_{l}(d+\delta)_{T}$ the $l$-order symbol of an operator $(d+\delta)_{T}$. An application of (2.1.4) in \cite{Wa1} shows that
\begin{align}
&\widetilde{{\rm Wres}}\Big[\pi^+\Big(c(u)c(v)c(w)(d+\delta)_{T_{2}}^{-p_1}\Big)\circ\pi^+(d+\delta)_{T_{2}}^{-p_2}\Big]\nonumber\\
=&\int_M\int_{|\xi|=1}{\rm Tr}_{\wedge^{*}(T^{*}M)}\Big[\sigma_{-n}( c(u)c(v)c(w)(d+\delta)_{T_{2}}^{-p_1-p_2})\Big]\sigma(\xi)\texttt{d}x+\int_{\partial
M}\Psi_{1},
\end{align}
and
\begin{align}
&\widetilde{{\rm Wres}}\Big[\pi^+\Big(c(u)\widehat{c}(v)\widehat{c}(w)(d+\delta)_{T_{2}}^{-p_1}\Big)\circ\pi^+(d+\delta)_{T}^{-p_2}\Big]\nonumber\\
=&\int_M\int_{|\xi|=1}{\rm Tr}_{\wedge^{*}(T^{*}M)}
\Big[\sigma_{-n}( c(u)\widehat{c}(v)\widehat{c}(w)(d+\delta)_{T}^{-p_1-p_2})\Big]\sigma(\xi)\texttt{d}x+\int_{\partial
M}\Psi_{2}.
\end{align}
and $r-k+|\alpha|+\ell-j-1=-n,r\leq-p_{1},\ell\leq-p_{2}$.

When $n=2m$, $ r-k-|\alpha|+l-j-1=-2m,~~r\leq-1,l\leq-2m+2$, and $r=l=-1,~k=|\alpha|=j=0,$ then we have the boundary terms:
\begin{align}
\Psi_{1}=&\int_{|\xi'|=1}\int^{+\infty}_{-\infty}
{\rm Tr}_{\wedge^{*}(T^{*}M)}
[ \sigma^+_{-1}(c(u)c(v)c(w)(d+\delta)_{T_{2}}^{-1})(x',0,\xi',\xi_n)\nonumber\\
&\times
\partial_{\xi_n}\sigma_{-2m+2}(d+\delta)_{T_{2}}^{-2m+2}(x',0,\xi',\xi_n)]\rm{d}\xi_{2m}\sigma(\xi')\rm{d}x',
\end{align}
and
\begin{align}
\Psi_{2}=&\int_{|\xi'|=1}\int^{+\infty}_{-\infty}
{\rm Tr}_{\wedge^{*}(T^{*}M)}
[ \sigma^+_{-1}(c(u)\widehat{c}(v)\widehat{c}(w)(d+\delta)_{T_{2}}^{-1})(x',0,\xi',\xi_n)\nonumber\\
&\times
\partial_{\xi_n}\sigma_{-2m+2}((d+\delta)_{T_{2}}^{-2m+2})(x',0,\xi',\xi_n)]\rm{d}\xi_{2m}\sigma(\xi')\rm{d}x'.
\end{align}
An easy calculation gives
  \begin{align}
 \pi^+_{\xi_n}\big( \sigma_{-1}(c(u)c(v)c(w)(d+\delta)_{T_{2}}^{-1})\big)=&
    c(u)c(v)c(w)\pi^+_{\xi_n}\big( \sigma_{-1} ((d+\delta)_{T_{2}}^{-1})\big)\nonumber\\
=&c(u)c(v)c(w) \frac{c(\xi')+ic(\rm{d}x_{n})}{2(\xi_n-i)} \nonumber\\
=& \frac{c(u)c(v)c(w)c(\xi')}{2(\xi_n-i)}+ \frac{ i c(u)c(v)c(w)c(\rm{d}x_{n})}{2(\xi_n-i)},
 \end{align}
where some basic facts and formulae about Boutet de Monvel's calculus which can be found  in Sec.2 in \cite{Wa1}.
By (2.16), we get
  \begin{align}
\partial_{\xi_n}\big(\sigma_{-2m+2}((d+\delta)_{T_{2}}^{-2m+2})\big)   (x_0)|_{|\xi'|=1}
=\partial_{\xi_n}\big(|\xi|^{2}\big)^{1-m}(x_0)=2(1-m)\xi_n(1+\xi_n^2)^{-m}.
\end{align}
By the relation of the Clifford action and $ {\rm{Tr}}(AB)= {\rm{Tr}}(BA) $, then
  \begin{align}
{\rm Tr}_{\wedge^{*}(T^{*}M)}[c(u)c(v)c(w)c({\rm{d}}x_{n})]=&\sum_{i,j,k=1}^{n}  u_{j}v_{k}w_{l}
{\rm Tr}_{\wedge^{*}(T^{*}M)}[c(e_{j})c(e_{k}) c(e_{l})c({\rm{d}}x_{n})]\nonumber\\
=&\sum_{i,j,k=1}^{n}  u_{j}v_{k}w_{l} (-\delta_{j}^{l}\delta_{k}^{n}+\delta_{j}^{n}\delta_{k}^{l}+\delta_{j}^{k}\delta_{n}^{l})
{\rm Tr}_{\wedge^{*}(T^{*}M)}[ {\rm Id}]\nonumber\\
=&\sum_{i,j,k=1}^{n}  (-u_{j}v_{n}w_{j}+u_{n}v_{k}w_{k}+w_{n}u_{j}v_{j})
{\rm Tr}_{\wedge^{*}(T^{*}M)}[ {\rm Id}]\nonumber\\
=&(u_{n}g(v,w)-v_{n}g(u,w)+w_{n}g(u,v)){\rm Tr}_{\wedge^{*}(T^{*}M)}[ {\rm Id}].
\end{align}
In the same way we have
  \begin{align}
&{\rm Tr}_{\wedge^{*}(T^{*}M)}[c(u)\widehat{c}(v)\widehat{c}(w)c({\rm{d}}x_{n})]\nonumber\\
=& {\rm Tr}_{\wedge^{*}(T^{*}M)}[c(u)\widehat{c}(w)\widehat{c}(v)c({\rm{d}}x_{n})]\nonumber\\
=&{\rm Tr}_{\wedge^{*}(T^{*}M)}[c(u)\big(-\widehat{c}(v)\widehat{c}(w)+2g(v,w)\big)c({\rm{d}}x_{n})]\nonumber\\
=&-{\rm Tr}_{\wedge^{*}(T^{*}M)}[c(u)\widehat{c}(v)\widehat{c}(w)c({\rm{d}}x_{n})]
+2g(v,w){\rm Tr}_{\wedge^{*}(T^{*}M)}[c(u)c({\rm{d}}x_{n})].
\end{align}
Then
  \begin{align}
{\rm Tr}_{\wedge^{*}(T^{*}M)}[c(u)\widehat{c}(v)\widehat{c}(w)c({\rm{d}}x_{n})]=&-u_{n}g(v,w){\rm Tr}_{\wedge^{*}(T^{*}M)}[ {\rm Id}].
\end{align}
From (5.7) and (5.8) we get
\begin{align}
& {\rm Tr}_{\wedge^{*}(T^{*}M)}
[ \sigma^+_{-1}(c(u)c(v)c(w)(d+\delta)_{T_{2}}^{-1})(x',0,\xi',\xi_n) \times
\partial_{\xi_n}\sigma_{-2m+2}((d+\delta)_{T_{2}}^{-2m+2})(x',0,\xi',\xi_n)] \nonumber\\
=& \frac{(1-m)\xi_n}{2(\xi_n-i)(1+\xi_n^2)^{m}}{\rm Tr}_{\wedge^{*}(T^{*}M)}[c(u)c(v)c(w)c(\xi')]\nonumber\\
&+ \frac{ i(1-m)\xi_n}{2(\xi_n-i)(1+\xi_n^2)^{m}}{\rm Tr}_{\wedge^{*}(T^{*}M)}[c(u)c(v)c(w)c(\rm{d}x_{n})]\nonumber\\
=& \frac{(1-m)\xi_n}{2(\xi_n-i)(1+\xi_n^2)^{m}}\sum_{i=1}^{n-1}\xi_i{\rm Tr}_{\wedge^{*}(T^{*}M)}[c(u)c(v)c(w)c(e_{i})]\nonumber\\
&+ \frac{ i(1-m)\xi_n}{2(\xi_n-i)(1+\xi_n^2)^{m}}(u_{n}g(v,w)-v_{n}g(u,w)+w_{n}g(u,v)){\rm Tr}_{\wedge^{*}(T^{*}M)}[ {\rm Id}].
\end{align}
From (5.7) and (5.9) we obtain
\begin{align}
& {\rm Tr}_{\wedge^{*}(T^{*}M)}
[ \sigma^+_{-1}(c(u)\widehat{c}(v)\widehat{c}(w)(d+\delta)_{T_{2}}^{-1}) (x',0,\xi',\xi_n) \times
\partial_{\xi_n}\sigma_{-2m+2}(d+\delta)_{T_{2}}^{-2m+2})(x',0,\xi',\xi_n)] \nonumber\\
=& \frac{(1-m)\xi_n}{2(\xi_n-i)(1+\xi_n^2)^{m}}{\rm Tr}_{\wedge^{*}(T^{*}M)}[c(u)\widehat{c}(v)\widehat{c}(w)c(\xi')]\nonumber\\
&+ \frac{ i(1-m)\xi_n}{2(\xi_n-i)(1+\xi_n^2)^{m}}{\rm Tr}_{\wedge^{*}(T^{*}M)}[c(u)\widehat{c}(v)\widehat{c}(w)c(\rm{d}x_{n})]\nonumber\\
=& \frac{(1-m)\xi_n}{2(\xi_n-i)(1+\xi_n^2)^{m}}\sum_{i=1}^{n-1}\xi_i{\rm Tr}_{\wedge^{*}(T^{*}M)}[c(u)c(v)c(w)c(e_{i})]\nonumber\\
&+ \frac{ i(1-m)\xi_n}{2(\xi_n-i)(1+\xi_n^2)^{m}}u_{n}g(v,w){\rm Tr}_{\wedge^{*}(T^{*}M)}[ {\rm Id}].
\end{align}
Substituting (5.11) into (5.4) we get
\begin{align}
\Psi_{1}=&\int_{|\xi'|=1}\int^{+\infty}_{-\infty}
{\rm Tr}_{\wedge^{*}(T^{*}M)}
[ \sigma^+_{-1}(c(u)c(v)c(w)(d+\delta)_{T_{2}}^{-1})(x',0,\xi',\xi_n)\nonumber\\
&\times
\partial_{\xi_n}\sigma_{-2m+2}((d+\delta)_{T_{2}}^{-2m+2})(x',0,\xi',\xi_n)]\rm{d}\xi_{2m}\sigma(\xi')\rm{ d}x'\nonumber\\
=&\int_{|\xi'|=1}\int^{+\infty}_{-\infty}\Big(
 \frac{(1-m)\xi_n}{2(\xi_n-i)(1+\xi_n^2)^{m}}\sum_{i=1}^{n-1}\xi_i{\rm Tr}_{\wedge^{*}(T^{*}M)}[c(u)c(v)c(w)c(e_{i})]\nonumber\\
&+ \frac{ i(1-m)\xi_n}{2(\xi_n-i)(1+\xi_n^2)^{m}}(u_{n}g(v,w)-v_{n}g(u,w)+w_{n}g(u,v)){\rm Tr}_{\wedge^{*}(T^{*}M)}[ {\rm Id}] \Big)
 \rm{d}\xi_{2m}\sigma(\xi')\rm{d}x'\nonumber\\
=& i(1-m)(u_{n}g(v,w)-v_{n}g(u,w)){\rm Tr}_{\wedge^{*}(T^{*}M)}[ {\rm Id}] {\rm{vol}}(S^{n-2})
 \int^{+\infty}_{-\infty}\frac{\xi_n}{2(\xi_n-i)(1+\xi_n^2)^{m}} \rm{d}\xi_{2m}\rm{d}x'\nonumber\\
 =& \frac{(2m-2)!(1-m)i2^{-2m+1}\pi}{m!(m-1)!}(u_{n}g(v,w)-v_{n}g(u,w)+w_{n}g(u,v)){\rm Tr}_{\wedge^{*}(T^{*}M)}[ {\rm Id}] {\rm{vol}}(S^{n-2})
{\rm d}\rm{vol}_{\partial_{M}}.
\end{align}
Substituting (5.12) into (5.5) we obtain
\begin{align}
\Psi_{2}=&\int_{|\xi'|=1}\int^{+\infty}_{-\infty}
{\rm Tr}_{\wedge^{*}(T^{*}M)}
[ \sigma^+_{-1}(c(u)\widehat{c}(v)\widehat{c}(w)(d+\delta)_{T_{2}}^{-1})(x',0,\xi',\xi_n)\nonumber\\
&\times
\partial_{\xi_n}\sigma_{-2m+2}((d+\delta)_{T_{2}}^{-2m+2})(x',0,\xi',\xi_n)]\rm{d}\xi_{2m}\sigma(\xi')\rm{ d}x'\nonumber\\
=& \frac{(2m-2)!(1-m)i2^{-2m+1}\pi}{m!(m-1)!}u_{n}g(v,w){\rm Tr}_{\wedge^{*}(T^{*}M)}[ {\rm Id}] {\rm{vol}}(S^{n-2})
{\rm d}\rm{vol}_{\partial_{M}}.
\end{align}
Summing up (5.13) and  Theorem 1.2 leads to the spectral torsion $\mathscr{\widetilde{T}}_{1}$ for manifold with boundary  as follows.
\begin{thm}
With the trilinear Clifford multiplication by functional of differential one-forms
$c(u),c(v),c(w)$, the spectral torsion for  de-Rham Hodge operator with torsion on manifold with boundary  equals to
\begin{align}
\mathscr{\widetilde{T}}_{1}(c(u),c(v),c(w))
=&\int_M (3-18m)2^{2m-1}T(u,v,w)
{\rm{vol}}(S^{2m-1}) {\rm d}\rm{vol}_{M}
+\int_{\partial_{M}}\frac{(2m-2)!(1-m)2^{-2m+1}\sqrt{-1}\pi}{m!(m-1)!}\nonumber\\
&\times(u_{n}g(v,w)-v_{n}g(u,w)+w_{n}g(u,v)){\rm Tr}_{\wedge^{*}(T^{*}M)}[ {\rm Id}] {\rm{vol}}(S^{2m-2})
{\rm d}\rm{vol}_{\partial_{M}}.
\end{align}
\end{thm}
Summing up (5.14) and  Theorem 1.3 leads to the spectral torsion $\mathscr{\widetilde{T}}_{2}$ for manifold with boundary  as follows.
\begin{thm}
With the trilinear Clifford multiplication by functional of differential one-forms
$c(u),\widehat{c}(v),\widehat{c}(w)$, the spectral torsion for  de-Rham Hodge operator with torsion  on manifold with boundary  equals to
\begin{align}
\mathscr{\widetilde{T}}_{2}(c(u),\widehat{c}(v),\widehat{c}(w))
=&\int_M (1-2m)2^{2m-1}T(u,v,w)
{\rm{vol}}(S^{2m-1}) {\rm d}\rm{vol}_{M}
+\int_{\partial_{M}}\frac{(2m-2)!(m-1)2^{-2m+1}\sqrt{-1}\pi}{m!(m-1)!}\nonumber\\
&\times u_{n}g(v,w){\rm Tr}_{\wedge^{*}(T^{*}M)}[ {\rm Id}] {\rm{vol}}(S^{2m-2})
{\rm d}\rm{vol}_{\partial_{M}}.
\end{align}
\end{thm}

\section{spectral forms for spectral triples}
Spectral triples $(\mathcal{A},\mathcal{H},D)$ are the noncommutative generalisations of compact spin manifolds
$M$. The aim of the noncommutative calculus is to adapt the notions of integration and differentiation
to spectral triples.
In \cite{DSZ2}, Dabrowski-Sitarz-Zalecki defined two spectral functionals for finitely summable spectral triples, which for the
canonical spectral triple over the spin manifold M allow to recover the metric and the Einstein
tensors, viewed as bilinear functionals over a pair of one-forms.
Similar to Section 2 and Section 3, we can construct the de-Rham Hodge Dirac triple
in this section.

In what follows, $\Gamma: E \rightarrow M$ is a real or complex vector bundle, linear means $\mathbb{R}$ or
$\mathbb{C}$ linear as appropriate, and $\Gamma(E)$ stands for the smooth sections of $E$ (as opposed to the
continuous sections).
For general spectral triples $(\mathcal{A},\mathcal{H},D)$, a connection on the space
of sections $\varepsilon\in\Gamma(E)$
will be defined as a linear map
$\nabla:\varepsilon\rightarrow \varepsilon\otimes_{\mathcal{A}}\Omega^{1}_{D}(\mathcal{A})$,
where $\Omega^{1}_{D}(\mathcal{A})$ consists of all finite linear combinations
of operators of the form $a[D,b]$.
The Leibniz rule then reads
$\nabla (\sigma\cdot f)=\sigma\otimes_{\mathcal{A}}df+ (\nabla\sigma)\cdot f$,
with $df=[D,f]$. As in the commutative case,
such connections always exist in multitude.

\begin{defn}\cite{DSZ2}
If $(\mathcal{A},\mathcal{H},D)$ is a $n$-summable spectral triple,let
$\Omega^{1}_{D}$ be the $\mathcal{A}$ bimodule of one forms generated by $\mathcal{A}$ and $[D, \mathcal{A}]$.
Moreover, assume there exists a
generalised algebra of pseudodifferential operators which contains $\mathcal{A}$, $D$,
$|D|^{l}$ for $l\in\mathbb{Z}$ with a
tracial state $Wres$  over this algebra (called a noncommutative residue),
 which identically vanishes on $T|D|^{-k}$ for any $k>n$ and a zero-order operator $T$
 (an operator in the algebra generated by $\mathcal{A}$ and
 $\Omega^{1}_{D}(\mathcal{A})$. Then for $u,v \in \Omega^{1}_{D}(\mathcal{A})$,
 the metric functional denoted by $\mathcal{G}_{D}(u,v):=Wres(uv|D|^{-n})$ and
  the Einstein functional denoted by $\widetilde{\mathcal{G}}_{D}(u,v):=Wres(u\{D,v\}D|D|^{-n})$.
\end{defn}
For de-Rham Hodge operator $d+\delta$, identifying for $dim M=n>2$,
 $\Omega^{1}_{d+\delta}(\mathcal{A})\simeq \Omega^{1}(M)$ via the correspondence of
  local components $u=\gamma^{p}u_{p}\leftrightarrow U=u_{p}e^{p} $ with respect to an orthonormal coframe
  $e^{p}$. Then Dabrowski-Sitarz-Zalecki \cite{DSZ} compute the metric functional $\mathcal{G}$.
 \begin{lem}\cite{DSZ}
for $u,v \in \Omega^{1}(M)$, the metric spectral functional reads
\begin{align}
\mathcal{G}=2^{n}w_{n-1}\int_{M}g(U,V)vol_{g}.
\end{align}
\end{lem}
\begin{lem}\cite{DSZ}
for $U,V \in \Omega^{1}(M)$, the Einstein functional reads
\begin{align}
\widetilde{\mathcal{G}}=\frac{2^{n}}{6}w_{n-1}\int_{M}G(U,V)vol_{g}.
\end{align}
\end{lem}
For de-Rham Hodge operator $d+\delta+\sqrt{-1}T_{1}$ and $d+\delta+\sqrt{-1}T_{3}$, for $dim M=n>2$,
 we have $\Omega^{1}_{d+\delta+\sqrt{-1}T_{1}}(\mathcal{A})\simeq \Omega^{1}(M)$
  and $\Omega^{1}_{d+\delta+\sqrt{-1}T_{3}}(\mathcal{A})\simeq \Omega^{1}(M)$.
  From Theorem 1.1, Lemma 6.2 and  Lemma 6.3, we obtain
\begin{thm}
With the bilinear Clifford multiplication by functional of differential one-forms
$c(u),c(v)$, the spectral triple for  de-Rham Hodge operator equals to
\begin{align}
{\rm Wres}\big(c(u)c(v)(d+\delta+\sqrt{-1}T_{1})^{-n+1}\big)
=\int_M (n-1)2^{n} \sqrt{-1}T(u,v )
\rm{vol}(S^{n-1}) {\rm d}\rm{vol}_{M},
\end{align}
where $T_{1}=\sum_{k<l }T(\widetilde{e_k},\widetilde{e_l})
c(\widetilde{e_j})c(\widetilde{e_l})$ is a two form.
\end{thm}
  From Theorem 1.4, Lemma 6.2 and  Lemma 6.3, we obtain
\begin{thm}
With the four-linear Clifford multiplication by functional of differential one-forms
$c(u),c(v)$, $ c(w),c(z)$, the spectral triple for de-Rham Hodge operator equals to
\begin{align}
{\rm Wres}\big(c(u)c(v)c(w)c(z)(d+\delta+\sqrt{-1}T_{3})^{-n+1}\big)
=\int_M \frac{(-n-1)}{3}2^{n-1}\sqrt{-1} T(u,v,w,z)
\rm{vol}(S^{n-1}) {\rm d}\rm{vol}_{M}.
\end{align}
where $T_{3}=\sum_{k<l <\alpha<\beta}T(\widetilde{e_k},\widetilde{e_l},\widetilde{e_\alpha},\widetilde{e_\beta})
c(\widetilde{e_k})c(\widetilde{e_l})c(\widetilde{e_\alpha})c(\widetilde{e_\beta})$ is a four form.
\end{thm}
From Definition 6.1, Theorem 6.4 and Theorem 6.5, we establish.
\begin{defn}
The spectral $l$-form functional for the spectral triple $(\mathcal{A},\mathcal{H},\widetilde{D})$
is a multilinear function for $u_{1},u_{2},\cdots,u_{l} \in\Omega_{D}^{l}(\mathcal{A})$ defined by
 \begin{align}
\mathscr{T}(u_{1},u_{2},\cdots,u_{l})=
{\rm Wres}\big(u_{1}u_{2}\cdots u_{l}\widetilde{D}^{-n+1}\big).
\end{align}
\end{defn}
Then the spectral triple is  the spectral $l$-form free, if $\mathscr{T}$
vanishes identically. By Section 5.2.1 \cite{DSZ2}, Dabrowski-Sitarz-Zalecki computed
 the metric and Einstein functionals for the conformally rescaled spectral triple on a
noncommutative tori.
\begin{exam}
Conformally  rescaled  Non-commutative tori, more examples see \cite{DSZ}.
\end{exam}
Then we have the following.
\begin{thm}\cite{DSZ}
The spectral triple    $(C^{\infty}(T_{n}^{\theta}),\mathcal{H},D_{k})$
for conformally rescaled Non-commutative tori is  the spectral $l$-form free.
\end{thm}

\section*{ AUTHOR DECLARATIONS}
Conflict of Interest:

The authors have no conflicts to disclose.
\section*{ Author Contributions}
Jian Wang: Investigation (equal); Writing- original draft (equal). Yong Wang: Investigation (equal); Writing-original draft (equal).
Mingyu Liu: Investigation (equal); Writing-original draft (equal).
\section*{DATA AVAILABILITY}
Data sharing is not applicable to this article as no new data were created or analyzed in this study.

\section*{ Acknowledgements}
The first author was supported by NSFC. 11501414.
The second author was supported by NSFC. 11771070. The authors also thank the referee for his (or her) careful reading and helpful comments.

\end{document}